\documentclass[reqno,oneside,12pt]{amsart}

\usepackage[T1]{fontenc}
\usepackage{times,mathptm}
\usepackage{amssymb,epsfig,verbatim,xypic}
\theoremstyle{plain}

\newtheorem{thm}{Theorem}[section]

\newtheorem{cor}[thm]{Corollary}
\newtheorem{pro}[thm]{Proposition}
\newtheorem{lem}[thm]{Lemma}
\newtheorem{proposition-principale}[thm]{Proposition principale}
\newtheorem{thm-principal}{Th\'eor\`eme principal}[section]

\theoremstyle{definition}
\newtheorem{defi}[thm]{Definition}

\newtheorem{eg}[thm]{Example}
\newtheorem{rem}[thm]{Remark}

\newenvironment{thm-A}
{\noindent{\bf Theorem A.}\it}{\\}

\newenvironment{thm-AA}
{\noindent{\bf Theorem A'.}\it}{\\}

\newenvironment{thm-B}
{\noindent{\bf Theorem B.}\it}{\\}

\newenvironment{thm-BB}
{\noindent{\bf Theorem B'.}\it}

\def\vv{\vspace{0.2cm}}

\def\C{\mathbf{C}}
\def\R{\mathbf{R}}
\def\Q{\mathbf{Q}}
\def\H{\mathbf{H}}
\def\Z{\mathbf{Z}}

\def\ii{{\sf{i}}}
\def\jj{{\sf{j}}}
\def\kk{{\sf{k}}}

\def\Nrd{{\text{Nrd}}}

\def\P{\mathbb{P}}

\def\Aut{{\sf{Aut}}}

\def\Int{{\sf{Int}}}

\def\PGL{{\sf{PGL}}\,}

\def\GL{{\sf{GL}}\,}
\def\SO{{\sf{SO}}\,}
\def\PSO{{\sf{PSO}}\,}

\def\U{{\sf{U}}\,}

\def\SL{{\sf{SL}}\,}
\def\Sp{{\sf{Sp}}\,}

\def\Mat{{\sf{Mat}}\,}

\def\Aff{{\sf{Aff}}\,}

\def\rk{{\sf{rk}}}

\def\aa{{\mathfrak{a}}}
\def\g{{\mathfrak{g}}}
\def\sll{{\mathfrak{sl}}}
\def\gll{{\mathfrak{gl}}}

\def\Ka{{\mathcal{K}}\,}
\def\Kabar{{\overline{\Ka}}}

\addtocounter{section}{0}             
\numberwithin{equation}{section}       

\begin{document}

\setlength{\baselineskip}{0.54cm}        
%
%
\title[Holomorphic actions and Zimmer Program]
{Holomorphic actions, Kummer examples, and Zimmer Program}
\date{October 2010}
\author{Serge Cantat and Abdelghani Zeghib}
\address{IRMAR (UMR 6625 du CNRS)\\ 
Universit{\'e} de Rennes 1 
\\ France}
\email{serge.cantat@univ-rennes1.fr}
\address{CNRS \\
UMPA \\
\'Ecole
Normale Sup\'erieure de Lyon\\
France}
\email{zeghib@umpa.ens-lyon.fr}

%
%

%
%

%
%

\begin{abstract} 

We classify compact K\"ahler manifolds $M$ of dimension $n\geq 3$ on which acts a lattice of an almost simple
real Lie group of rank $\geq n-1$. This provides a new line in the so-called Zimmer program, and characterizes
certain type of complex tori by a property of their automorphisms groups. 

\vspace{0.1cm}

\noindent{\sc{R\'esum\'e.}} Nous classons les vari\'et\'es complexes compactes k\"ahl\'eriennes $M$ de dimension $n\geq 3$
munies d'une action d'un r\'eseau $\Gamma$ dans un groupe de Lie r\'eel presque simple de rang $n-1$. Ceci compl\`ete
 le programme de Zimmer dans ce cadre, et caract\'erise certains tores complexes compacts par des
propri\'et\'es de leur groupe d'automorphismes.

\end{abstract}

\maketitle

\setcounter{tocdepth}{1}
\tableofcontents
%
%

\section{Introduction}
%
%

\subsection{Zimmer Program}
Let $G$ be an almost simple real Lie group. The {\bf{real 
rank}} $\rk_\R(G)$ of $G$ is the
dimension of a maximal abelian subgroup $A$ of $G$ that acts by $\R$-diagonalizable endomorphisms in the adjoint
representation of $G$ on its Lie algebra $\g$.   
When $\rk_\R(G)$ is at least $2$, we shall say
that $G$ is a {\bf{higher rank}} almost simple Lie group.
Let $\Gamma$ be a {\bf{lattice}} in $G$; by
definition, $\Gamma$ is a discrete subgroup of $G$ such that
$G/\Gamma$ has finite Haar volume. 
Margulis superrigidity theorem
implies that all finite dimensional linear representations of $\Gamma$ are built from 
representations in unitary groups and representations of the Lie group $G$
itself. In particular, there is no faithful linear representation of $\Gamma$
in dimension $\leq \rk_\R(G)$ (see \cite{Margulis:book}).

Zimmer's program predicts that a similar picture should hold
for actions of $\Gamma$ by diffeomorphims on compact manifolds, 
at least when the dimension $\dim(V)$ of the manifold $V$ is close to the minimal 
dimension of non trivial linear representations of $G$ (see \cite{Fisher:preprint}).   
For instance, a central conjecture predicts that lattices
in simple Lie groups of rank $n$ do not act faithfully on compact manifolds of dimension less than
$n$ (see \cite{Zimmer:ICM, Zimmer:1984, Zimmer:Mostow, Ghys:1999}).

In this article, we pursue the study of Zimmer's program in the holomorphic, k\"ahlerian,
setting, as initiated in \cite{Cantat:ENS} and \cite{Cantat-Zeghib}. 

\subsection{Automorphisms} Let $M$ be a compact complex manifold of dimension~$n$.
By definition, diffeomorphisms of $M$ which are holomorphic are called {\bf{automorphisms}}.   
According  to Bochner and Montgomery \cite{Bochner-Montgomery:1946, Campana-Peternell:survey}, 
the group $\Aut(M)$ of all automorphisms of $M$
is a complex Lie group, the Lie algebra
of which is the algebra of holomorphic vector fields on $M$. 
Let $\Aut(M)^0$ be the connected component
of the identity in $\Aut(M)$, and 
\[
\Aut(M)^\sharp=\Aut(M)/\Aut(M)^0
\]
be the group of connected components. This group can be infinite, and is hard to 
describe: For example, it is not known whether there exists a  compact complex
manifold $M$ for which $\Aut(M)^\sharp$ is not finitely generated. 

When $M$ is a K\"ahler manifold, Lieberman and Fujiki proved that $\Aut(M)^0$ has finite index in the kernel of the action of $\Aut(M)$
on the cohomology of $M$ (see \cite{Fujiki:1978, Lieberman:1978}). Thus, if a  subgroup $\Gamma$ of $\Aut(M)$ 
embeds into  $\Aut(M)^\sharp$, the action of $\Gamma$  on the 
cohomology of $M$ has finite kernel; in particular, the group $\Aut(M)^\sharp$ almost embeds in the group ${\text{Mod}}(M)$
of isotopy classes of smooth diffeomorphisms of $M$. 
When $M$ is simply connected, ${\text{Mod}}(M)$
is naturally described as the group of integer matrices in a linear algebraic group~(\cite{Sullivan:1977}). 
Thus, $\Aut(M)^\sharp$  sits naturally in an arithmetic lattice. 
Our main result goes in the other direction: it describe the largest possible lattices contained in $\Aut(M)^\sharp$.

\subsection{Rigidity and Kummer examples} The main example which provides large groups $\Gamma\subset \Aut(M)^\sharp$ is 
given by linear actions on tori, and on quotient of tori (see \cite{Cantat-Zeghib}, \S 1.2). For instance, if
$\Lambda_0$ is a lattice in $\C$, the group $\SL_n(\Z)$ acts on the torus  $A=(\C/\Lambda_0)^n$; since this
action commutes with multiplication by~$-1$, $\SL_n(\Z)$ also acts on the quotient $M_0=A/\langle -1\rangle$ and on the
smooth $n$-fold $M$ obtained by blowing up the $4^n$ singularities of $M_0$. The following 
definition, which is taken from \cite{Cantat:Compositio, Cantat-Zeghib}, provides a common denomination for all these examples.

\begin{defi}
{\sl{Let $\Gamma$ be a group, and $\rho:\Gamma\to \Aut(M)$ a morphism
into the group of automorphisms of a compact complex manifold $M.$ This
morphism is a {\bf{Kummer example}} (or, equivalently, is of {\bf{Kummer type}}) if there
exists
\begin{itemize}
\item a birational morphism $\pi:M\to M_0$ onto  an orbifold $M_0$, 
\item a finite orbifold cover $\epsilon:A\to M_0$ of $M_0$ by a torus $A,$ and
\item a morphism $\eta:\Gamma \to \Aut(A)$
\end{itemize} 
such that 
$
\epsilon \circ \eta(\gamma)
=
(\pi\circ \rho(\gamma)\circ \pi^{-1})\circ \epsilon  
$
for all $\gamma$ in $\Gamma.$}}
\end{defi}

The notion of {\bf{orbifold}} used  in  this text refers to 
compact complex analytic spaces with a finite number of singularities of
quotient type; in other words, $M_0$ is locally the quotient of $(\C^n, 0)$
by a finite group of linear tranformations (see Section \ref{par:Orbifolds}). 

Since automorphisms of a torus $\C^n/\Lambda$ are covered by affine transformations of $\C^n$,
all Kummer examples are covered by the action of affine transformations on 
the affine space.  

The following statement is our main theorem. 
It confirms Zimmer's program, in its strongest versions, for holomorphic actions on compact K\"ahler manifolds:
 We get a precise description of all possible actions 
of lattices $\Gamma\subset G$ for $\rk_\R(G)=\dim_\C(M)$ but also for $\rk_\R(G)=\dim_\C(M)-1$.

\vspace{0.209cm}

\noindent{\bf{Main Theorem.}}{\it{
Let $G$ be an almost simple real Lie group and $\Gamma$ be a lattice in $G$. 
Let $M$ be a compact K\"ahler manifold of dimension $n\geq 3$. 
Let $\rho:\Gamma\to \Aut(M)$ be an injective morphism. 
Then, the real rank $\rk_\R(G)$ is at most equal to the complex dimension of $M$.
\begin{itemize}
\item[(1)] If $\rk_\R(G)=\dim(M)$, then $G$ is locally isomorphic to $\SL_{n+1}(\R)$
or $\SL_{n+1}(\C)$ and $M$ is biholomorphic to the  projective space~$\P^n(\C)$.
\item[(2)] If $\rk_\R(G)=\dim(M)-1$, there exists a finite index subgroup $\Gamma_0$
in $ \Gamma$ such that either

{\rm{(2-a)}}  $\rho(\Gamma_0)$ is contained in $\Aut(M) ^0$, or

{\rm{(2-b)}}  $G$ is locally isomorphic to $\SL_{n}(\R)$ or $\SL_{n}(\C)$, and
the morphism $\rho: \Gamma_0 \to \Aut(M)$ is a Kummer example.
\end{itemize}
}}

\vspace{0.1cm}

Moreover, all examples corresponding to assertion (2-a) are described in Section \ref{par:list2b}
and all  Kummer examples of assertion (2-b) are described in Section \ref{par:Kummer}.
In particular, for these Kummer examples,  the complex torus $A$ associated to $M$ and the lattice 
$\Gamma$ fall in one of the following three possible examples:

\vspace{0.1cm}

$\bullet$  $\Gamma\subset \SL_n(\R)$ is commensurable
to $\SL_n(\Z)$ and  $A$ is isogeneous to the product of $n$ copies of an elliptic curve $\C/\Lambda$;

\vspace{0.1cm}

$\bullet$  $\Gamma\subset\SL_n(\C)$ is commensurable to $\SL_n({\mathcal{O}}_d)$ where ${\mathcal{O}}_d$
is the ring of integers in $\Q(\sqrt{d})$ for some negative integer $d$, and $A$ is isogeneous to the product of $n$ copies of the elliptic curve $\C/{\mathcal{O}}_d$;

\vspace{0.1cm}

$\bullet$  In the third example, $n=2k$ is even.  There are integers $a$ and $b$ such that $A$ is isogeneous to 
the product of $k$ copies of the abelian surface $\C^2/\H_{a,b}(\Z)$, where $\H_{a,b}$ is the division algebra of quaternions
over the rational numbers $\Q$ defined by the basis $(1,\ii,\jj,\kk)$, with 
\[
\ii^2=a, \, \jj^2=b, \, \ii\jj=\kk=-\jj\ii.
\]
Moreover, the group $\Gamma$ is a lattice in $\SL_n(\R)$ commensurable to the group of automorphisms
of the abelian group $\H_{a,b}(\Z)^k$ that commute to the diagonal action of $\H_{a,b}(\Z)$ by left multiplications (see \S \ref{par:Kummer}).

\vspace{0.1cm}

As a consequence, $\Gamma$ is not cocompact, $A$ is an abelian variety and $M$ is projective. 
This theorem extends
the main result of \cite{Cantat-Zeghib} from dimension $3$ to all dimensions $n\geq 3$ when $G$ is almost simple; the strategy 
is different, more concise, but slightly less precise. 

\subsection{Strategy of the proof and complements} 
After a few preliminary facts (\S \ref{PART:2}), the proof of the Main Theorem starts in \S \ref{PART:3}: 
Assertion (1) is proved, and a complete list of all possible pairs $(M,G)$ in assertion (2-a) is obtained. 
This makes use of a previous result on  Zimmer conjectures in the holomorphic setting 
(see \cite{Cantat:ENS}), and classification of homogeneous or quasi-homogeneous spaces
(see \cite{Akhiezer:book, Gilligan-Huckleberry:1981, Huckleberry-Snow:1982}). On our way, we describe
$\Gamma$-invariant analytic subsets $Y\subset M$.

The core of the paper proves that assertion (2-b) is satisfied when  the image $\rho(\Gamma_0)$ is not contained
in $\Aut(M)^0$ and $\rk_\R(G)=\dim(M)-1$. 

In that case, $\Gamma$ acts almost faithfully on the cohomology of $M$, and this linear representation extends
to a continuous  representation of $G$ on $H^*(M,\R)$. Section \ref{PART:4} shows that $G$ preserves a non-trivial cone contained in the
closure of the K\"ahler cone $\Ka(M)\subset H^{1,1}(M,\R)$; this general fact holds for all linear representations of semi-simple Lie
groups $G$ for which a lattice $\Gamma\subset G$ preserves a salient cone. Section \ref{PART:4} can be skipped in a first reading.

Then, in \S \ref{PART:5}, we apply ideas of Dinh, Sibony and Zhang together with 
representation theory. We fix a Cartan subgroup $A$ in $G$ and study the eigenvectors
of $A$ in the $G$-invariant cone: Hodge index theorem constrains the set of weights
and eigenvectors. When there is no $\Gamma$-invariant analytic subset of positive
dimension, Yau's Theorem can then be used to prove that $M$ is a torus. 

To conclude the proof, we then show that invariant analytic
subsets can be blown down to quotient singularities, and we apply Hodge and Yau's theorems in the orbifold setting. 
This makes use of Section \ref{PART:3}. 

Section \ref{par:Kummer} lists all tori of dimension $n$ with an action of a lattice in 
a simple Lie group of rank $n-1$.  
Since Sections \ref{par:list2b} and \ref{par:Kummer} provide complements to the Main Theorem, we recommend to skip them 
in a first reading.

\subsection{Aknowledgment} Thanks to Michel Brion, Jean-Pierre Demailly, Igor Dolgachev, St\'ephane Druel, Jean-Fran\c{c}ois Quint for nice
discussions, comments, and ideas. Demailly provided the proof of Theorem \ref{thm:DP} while  Brion and Dolgachev helped us clarify Section \ref{part:homorank}.

%
%

\section{Cohomology, Hodge theory, Margulis extension}\label{PART:2}
%
%

Let $M$ be a connected, compact, K\"ahler manifold of complex dimension $n$.

\subsection{Hodge Theory and cohomological automorphisms}

\subsubsection{Hodge decomposition} Hodge theory implies that 
the cohomology groups $H^k(M,\C)$ decompose into direct sums
\[
H^k(M,\C)= \bigoplus_{p+q=k} H^{p,q}(M,\C), 
\]
where cohomology classes in $H^{p,q}(M,\C)$ are represented by closed
forms of type $(p,q).$ This bigraded structure is compatible with the cup product. 
Complex conjugation permutes $H^{p,q}(M,\C)$
with $H^{q,p}(M,\C).$ In particular, the cohomology groups $H^{p,p}(M,\C)$ 
admit a real structure, the real part of which is 
\[
H^{p,p}(M,\R)= H^{p,p}(M,\C)\cap H^{2p}(M,\R).
\]
If $[\kappa]$ is a K\"ahler class (i.e. the cohomology class of a K\"ahler form),
then $[\kappa]^p \in H^{p,p}(M,\R)$ for all $p.$

\subsubsection{Notation}

In what follows, the vector space $H^{1,1}(M,\R)$ is denoted~$W$.

\subsubsection{Primitive classes and Hodge index theorem}

Let $[\kappa]\in W$ be a {\bf{K\"ahler class}}, i.e. the class of a K\"ahler form (alternatively, K\"ahler classes are also called {\bf{ample}} classes).
The set of primitive classes with respect to $[\kappa]$
is the vector space of classes $[u]$ in $W$ such that 
\[
\int_M[\kappa]^{n-1}\wedge [u] = 0.
\]
Hodge index theorem implies that the quadratic form 
\[
([u],[v])\mapsto \int_M [\kappa]^{n-2}\wedge [u]\wedge [v]
\]
is negative definite on the space of primitive forms (see \cite{Voisin:book}, \S 6.3.2).
We refer the reader to~\cite{Dinh-Sibony:Duke}, \cite{Dinh-Nguyen:2009} and~\cite{Zhang:Invent} for stronger results
and consequences on groups of automorphisms of $M$.

\subsubsection{Cohomological automorphisms}

\begin{defi}
A {\bf{cohomological automorphism}} of $M$ is a  linear isomorphism of the real vector space $H^*(M,\R)$ that
preserves the Hodge decomposition, the cup product, and the Poincar\'e duality. 
\end{defi}

Note that cohomological automorphisms are not assumed to preserve the set of K\"ahler classes or
the lattice $H^*(M,\Z),$ as automorphisms $f^*$ with $f\in \Aut(M)$ do.

\subsection{Nef cone and big classes}\label{par:NefBig}

Recall that a convex cone  in a real vector space 
is {\bf{salient}} when it does not contain any line: In other words, a salient cone is
strictly contained in a half space.

The {\bf{K\"ahler cone}} of $M$ is the subset $\Ka(M)\subset W$ of K\"ahler classes.
This set is an open convex cone;  its closure
${\overline{\Ka}}(M)$ is a strict and closed convex cone, the interior of which coincides
with $\Ka(M).$ We shall say that
${\overline{\Ka}}(M)$ is the cone of {\bf{nef}} cohomology classes of type $(1,1)$. 
All these cones are invariant under the action of $\Aut(M).$  

A class $[\omega]$ in $H^{1,1}(M,\R)$ is {\bf{big and nef}} if it is nef and $\int_M \omega^n>0$. 
The cone of big and nef classes plays an important role in this paper.  

\begin{thm}[Demailly and Paun]\label{thm:DP}
Let $M$ be a compact K\"ahler manifold, and $[\omega]\in H^{1,1}(M,\R)$ be a big and nef class
which is not a K\"ahler class. Then 
\begin{enumerate}
\item there exists an irreducible analytic subset $Y\subset M$ of positive dimension  such that
\[
\int_Y\omega^{\dim(Y)}=0\, ;
\]
\item the union of all these analytic subsets $Y$ is a proper Zariski closed subset $Z\subset M$. 
\end{enumerate}
\end{thm}

\begin{proof}
The existence of $Y$ in property (1) follows from Theorem 0.1 in \cite{Demailly-Paun:2004}.

Let us now prove property (2).
Theorem 0.5 in \cite{Demailly-Paun:2004} shows that the class $[\omega]$ is represented by a closed
positive current $T$ which is smooth in the complement of a proper analytic subset
$Z\subset M$, has logarithmic poles along $Z$, and is bounded from below by a K\"ahler form,
i.e. $T\geq \kappa$ for some K\"ahler form $\kappa$ on~$M$.
Our goal is to show that all irreducible analytic subsets $Y\subset M$ of positive dimension
that satisfy property (1) are contained in $Z$. We assume that $Y$ is not contained in $Z$
and $\dim(Y)>0$, and we want to show that the integral of $[\omega]^{\dim(Y)}$ on $Y$ is positive.
In order to compute this integral, we represent $\omega$ by $T$ and regularize
$T$ in order to take its $\dim(Y)$-exterior power. 

Let $\alpha$ be a smooth and closed form of type $(1,1)$ which represents the class $[\omega]$.
Let $C>0 $ be a constant such that $\alpha \geq -C \kappa$. Write $T$ as
\[
T=\alpha+ \frac{i}{\pi}\partial \overline\partial \, \psi \geq  \kappa,
\]
and consider the sequence of truncated currents $T_a$, $a>0$, defined by 
\[
T_a=\alpha + \frac{i}{\pi}\partial \overline\partial \, \max(\psi, -a). 
\]
On the set $\psi>-a$, $T_a$ coincides with $T$ and thus $T_a\geq \kappa$; on the set $\psi<-a$ it coincides
with $\alpha$.  In particular, $T_a\geq  -C\kappa$ on $M$.
Since $\psi$ has logarithmic poles along $Z$, the sets $\psi<-a$ are
contained in smaller and smaller neighborhoods of $Z$ when $a$ goes to $\infty$.

Since $\psi$ is locally  the difference of a smooth function and a plurisubharmonic 
function, $\psi$ is upper semi-continuous and, as such, is bounded from above. Thus, 
$T_a$ has bounded local potentials, and its Monge-Amp\`ere products can be computed
on any analytic subset of $M$ by Bedford-Taylor technique (see \cite{Bedford-Taylor:1976}). 

Since the cohomology class of $T_a$ is equal to the class of $T$ we have 
\[
\int_Y [\omega]^{\dim(Y)} = \int_Y T_a^{\dim(Y)} \geq  \int_{Y\cap\{\psi>-a+1\} } \kappa^{\dim(Y)} - \int_{Y\cap\{\psi<-a + 1\}} (C\kappa)^{\dim(Y)}.
\]
The first term of the right hand side of this inequality goes to zero when $a$ goes to $-\infty$. The second term converges to the volume
of $Y$ with respect to $\kappa$. This concludes the proof.\end{proof}

\subsection{Margulis rigidity and extension}

\vv

Let $H$ be a group. A property is said to hold {\bf{virtually}} for $H$ if
a finite index subgroup of $H$ satisfies this property. Similarly, a morphism
$h:\Gamma\to L$ from a subgroup $\Gamma$ of $H$ to a group $L$ 
virtually extends to $H$ if there is a finite index subgroup $\Gamma_0$ in $\Gamma$
and a morphism ${\hat{h}}:H\to L$ such that $\hat{h}$ coincides with $h$ on the subgroup $\Gamma_0$.

The following theorem is one version of the superrigidity phenomenum for linear representations of lattices (see \cite{Margulis:book} or \cite{VGS:EMS}). 

\begin{thm}[Margulis]\label{thm:SuperMargu}
Let $G$ be a semi-simple connected Lie group with finite center, with rank at least $2,$ and without non trivial compact factor.
Let $\Gamma\subset G$ be an irreducible lattice. 
Let $h:\Gamma \to \GL_k(\R)$ be a linear representation of $\Gamma.$
The Zariski closure of $h(\Gamma)$ is a semi-simple Lie group; if this
Lie group does not have any infinite compact factor, then $h$ virtually extends to a (continuous) linear representation 
${\hat{h}}:G\to \GL_k(\R)$.\end{thm}

Another important statement due to Margulis asserts that irreducible lattices $\Gamma$ in higher rank semi-simple Lie groups
are "almost simple": If $\Gamma'$ is a normal subgroup of $\Gamma$, either $\Gamma'$ is finite or $\Gamma'$ has
finite index in $\Gamma$. Thus, if $\rho$ is a morphism from $\Gamma$ to a group $L$ with infinite image, then $\rho$ is virtually faithful 
(see \cite{Margulis:book} or \cite{VGS:EMS}).

As explained in \cite{Cantat-Zeghib}, Margulis theorems, Lieberman-Fujiki theorem, 
and the fact that the action of $\Aut(M)$ on $H^*(M,\R)$ preserves the lattice $H^*(M,\Z)$ 
imply the following proposition. 

\begin{pro}\label{pro:ExtensionCohomology}
Let  $G$ and $\Gamma$ be as in theorem \ref{thm:SuperMargu}.
Let $\rho:\Gamma \to \Aut(M)$ be a representation into the
group of automorphisms
of a compact K\"ahler manifold $M.$ Let $\rho^* :\Gamma\to \GL(H^*(M,\Z))$ be the
induced action on the cohomology ring of $M.$ 
\begin{itemize}
\item[(a)] If the image of $\rho^*$ is
infinite, then $\rho^*$ virtually extends to a representation 
${\hat{\rho^*}}:G\to \GL(H^*(M,\R))$ by cohomological automorphisms.

\item[(b)] If the image of $\rho^*$ is finite, the image of $\rho$ is virtually contained
in $\Aut(M)^0.$ 
\end{itemize}
\end{pro}

\subsection{Orbifolds}\label{par:Orbifolds}

In this paper, an orbifold $M_0$ of dimension $n$ is a compact complex analytic space 
with a finite number of quotient singularities $q_i$; in a neighborhood of each $q_i$, $M_0$ is
locally isomorphic to the quotient of $\C^n$ near the origin by a finite group of linear
transformations. All examples of orbifolds considered in this paper are locally isomorphic
to $\C^n /\eta_i$ where $\eta_i$ is a scalar multiplication of finite order $k_i$.  Thus, the
singualrity $q_i$ can be resolved by one blow-up: The point $q_i$ is then replaced by a hypersurface
$Z_i$ which is isomorphic to $\P^{n-1}(\C)$ with normal bundle ${\mathcal{O}}(-k_i)$. 

All classical objects from complex differential geometry are defined on $M_0$ 
as follows. Usual definitions are applied on the smooth part $M_0\setminus\{ q_1, ..., q_k\}$ and,  around
each singularity $q_i$, one requires that the objects come locally from  $\eta_i$-invariant 
objects on $\C^n$. 
Classical facts, like Hodge decomposition, Hodge index theorem, Yau theorem, remain
valid in the context of orbifolds. The reader will find more details  in~\cite{Campana:2004, Zaffran-Wang:2009}.

%
%

\section{Lie group actions and invariant analytic subsets}\label{PART:3}
%
%

\subsection{Homogeneous manifolds}\label{part:homorank}

The following theorem is a direct consequence of the classification of maximal subgroups in simple Lie groups (see
\cite{Onishchik-Vinberg:EMS3}, chapter 6, or Section \ref{par:list2b} below).

\begin{thm}\label{thm:homogeneous}
Let $H$ be a connected almost simple complex Lie group of rank $\rk_\C (H) = n$. If $H$ acts faithfully  and holomorphically
on a connected compact complex manifold $M$ of dimension $\leq n$ then, up to holomorphic conjugacy, 
$M$ is the projective space $\P^n(\C)$, $H$ is locally isomorphic to $\PGL_{n+1}(\C)$, and the action of $H$ on
$M$ is the standard action by linear projective transformations.
\end{thm}

Following a suggestion by  Brion and Dolgachev, we sketch a proof that does not use the 
classification of maximal subgroups of Lie groups.
Let $A$ be a Cartan subgroup in $H$. Since $H$ 
has rank $n$ this group is isomorphic to the multiplicative group
$(\C^*)^n$. The action of $A$ on $M$ is faithful; this easily implies that $\dim(M)=n$ and that 
$A$ has an open orbit. 
Thus $M$ is a toric variety of dimension $n$ with respect to the action of the multiplicative group~$A$.
In particular, there is no faithful action of $H$ on compact complex manifolds of dimension less than~$n$. 
Since $H$ is almost simple and connected, all actions of $H$ in dimension $<n$ are trivial. 

As a corollary, $H$ acts transitively on $M$, because otherwise $H$ has a proper Zariski closed orbit: This orbit
has dimension $< n$ and, as such, must be a point $m\in M$; the action of $H$ at $m$ can be linearized, 
and gives a non-trivial morphism from $H$ to $\GL(T_mM)\simeq \GL_n(\C)$, in contradiction with $\rk_\C(H)=n$. Thus $M=H/L$ for some closed subgroup $L$. 

Since $H/L$ is compact, $L$ is contained in a parabolic subgroup $P$ (see \cite{Akhiezer:book}). Since the dimension
of $M$ is the smallest positive dimension of a $H$-homoge\-neous space, $P=L$ and $P$ is a maximal 
parabolic subgroup.
Since $P$ is maximal, the Picard number of $M$ is equal to $1$ (see \cite{Akhiezer:book}, \S 4.2). 

As a consequence, $M$ is a smooth toric variety with Picard number $1$ and, as such, is isomorphic to $\P^n(\C)$ (see \cite{Fulton:book}). 
Since the group of automorphisms of $\P^n(\C)$ is the rank $n$ group  $\PGL_{n+1}(\C)$, the conclusion 
follows.

\subsection{First part of the Main Theorem}\label{par:firstpart}

Let us apply Theorem \ref{thm:homogeneous}. Let $\Gamma$ be a lattice in an almost simple real Lie group $G$. 
Assume that $\Gamma$ acts faithfully on a connected compact K\"ahler manifold $M$, with 
$\dim_\C(M)\leq \rk_\R(G)$. By \cite{Cantat:ENS}, the dimension of $M$ is equal to the rank of $G$ and
the image of $\Gamma$ in $\Aut(M)$ is virtually contained in $\Aut(M)^0$.
 Hence, we can assume that
the action of $\Gamma$ on $M$ is given by an injective morphism $\rho\colon \Gamma\to \Aut(M)^0$. 
As explained in \cite{Cantat:ENS}, the complex Lie group $\Aut(M)^0$ contains a copy of an almost simple
complex Lie group $H$ with $\rk_\C(H)\geq \rk_\R(G)$. More precisely, if $\rho(\Gamma)$
is not relatively compact in $\Aut(M)^0$, one apply Theorem \ref{thm:SuperMargu} to extend the morphism $\rho$
virtually to a morphism $\hat{\rho}\colon G\to \Aut(M)^0$; if the image of $\rho$ is relatively compact, then
another representation $\rho'\colon \Gamma\to \Aut(M)^0$ extends virtually to $G$; in both cases, the Lie algebra of $H$ 
is the smallest complex Lie subalgebra containing ${\text{d}}\hat{\rho}_{Id}({\mathfrak{g}})$.

Theorem \ref{thm:homogeneous} shows
that $M$ is the projective space $\P^n(\C)$ and   $\Aut(M)$ coincides with $\PGL_{n+1}(\C)$
(and thus with $H$). As a consequence, the group $G$ itself is locally isomorphic to $\SL_{n+1}(\R)$
or $\SL_{n+1}(\C)$.

Summing up, the inequality $\dim_\C(M)\geq \rk_\R(G)$ as well as property (1) in the Main Theorem
have been proved.

\subsection{Invariant analytic subsets}

Let us now study $\Gamma$-invariant analytic subsets $Z\subset M$ under the assumption of assertion (2)
 in the Main Theorem; in particular $\dim_\C(M)=\rk_\R(G)+1$.
Let $Z$ be a $\Gamma$-invariant complex analytic subset. Assume, first, that 
(i) $Z$ is irreducible and (ii) $Z$ has positive dimension.

\subsubsection{Singularities} 

If $\dim(Z)<n-1$, part (1) of the Main Theorem implies that a finite index subgroup of $\Gamma$ 
fixes $Z$ pointwise. If the set $Z$ is not smooth, its singular locus is $\Gamma$-invariant and has dimension $\leq n-2$.
Hence, changing $\Gamma$ into a finite index subgroup, we assume that $\Gamma$ fixes the singular locus
of $Z$ pointwise as well as $Z$ itself if its codimension is larger than $1$.
 
If $\Gamma$ fixes a point $q\in Z$, the image of the morphism $\delta_q\colon \Gamma\to \GL(T_qM)$
defined by the differential at $q$, i.e. by
\[
\delta_q(\gamma)={\text{d}}\gamma_{q},
\]
preserves the tangent cone of $Z$ at $q$; in particular, the Zariski closure of $\delta_q(\Gamma)$ in 
$\PGL(T_q M)$ is a proper algebraic subgroup of $\PGL(T_q M)$.
Since proper algebraic subgroups of $\PGL_n(\C)$
have rank less than $n-1=\rk_\R(G)$, Margulis rigidity theorem implies that the image of $\delta_q$ is finite.
These facts provide the following alternative: 
\begin{itemize}
\item Either $Z$ is a smooth hypersurface,
\item or $Z$ contains a fixed point $q$ for which the morphism $\delta_q$ has finite image.
\end{itemize}
From \cite{Cairns-Ghys:1997}, the action of $\Gamma$ in a 
neighborhood of  a fixed point $q$ can be linearized. Thus, in  the second alternative, a finite index subgroup of $\Gamma$ acts trivially in a neighborhood of $q$. 
Since  the action of $\Gamma$ is holomorphic and $M$ is connected, this contradicts the faithfulness of
the morphism $\Gamma\to \Aut(M)$. We deduce that {\sl{all irreducible $\Gamma$-invariant analytic subsets of positive dimension
are smooth hypersurfaces}}.

\subsubsection{Geometry of $Z$}

By the first part of the Main Theorem (see \S \ref{par:firstpart}), a smooth invariant hypersurface is a copy of $\P^{n-1}(\C)$ on which $\Gamma$
acts as a Zariski dense subgroup of $\PGL_{n}(\C)$. Such a subgroup does not preserve any non empty algebraic
subset. Thus, $Z$ does not intersect any other irreducible $\Gamma$-invariant subset. Replacing $\Gamma$
by finite index subgroups, one can now apply this discussion to all $\Gamma$-invariant analytic subsets:  

\begin{pro}\label{pro:ias}
Let $\Gamma$ be a lattice in an almost  simple Lie group of rank $n-1\geq 2$. 
If $\Gamma$ acts faithfully by holomorphic transformations on a compact complex
manifold $M$ of dimension $n$, any $\Gamma$-invariant analytic subset $Z\subset M$
is a disjoint union of isolated points and smooth hypersurfaces isomorphic to $\P^{n-1}(\C)$.
\end{pro}

\subsubsection{Contraction of $Z$}

Section 3.2 in \cite{Cantat-Zeghib} can now be applied almost word by word to show the following result. 

\begin{thm}\label{thm:ih}
Let $\Gamma$ be a lattice in an almost simple Lie group $G$. Assume that 
$\Gamma$ acts faithfully on a connected compact K\"ahler manifold $M$,  
\[
\rk_\R(G)=\dim_\C(M)-1,
\]
and the image of $\Gamma$ in $\Aut(M)$  is not virtually contained in $\Aut(M)^0$.

Let $Z$ be the union of all $\Gamma$-invariant analytic subsets $Y\subset M$ with
positive dimension. Then $Z$ is the union of a finite number of disjoint copies of the projective space
$Z_i=\P^{n-1}(\C)$. Moreover there exists a birational morphism $\pi\colon M\to M_0$
onto a compact K\"ahler orbifold $M_0$ such that 
\begin{itemize}
\item[(1)] $\pi$ contracts all $Z_i$ to points $q_i\in M_0$;
\item[(2)] around each point $q_i$, the orbifold $M_0$ is either smooth, or
locally isomorphic to a quotient of $(\C^n, 0)$ by a finite order scalar multiplication;
\item[(3)] $\pi$ is an isomorphism from the complement of $Z$ to the complement of the points $q_i$;
\item[(4)] $\pi$ is equivariant: The group $\Gamma$ acts on $M_0$ in such a way that
$\pi\circ \gamma=\gamma\circ \pi$ for all $\gamma$ in $\Gamma$.
\end{itemize}
\end{thm}

\subsection{Lie group actions in case $\rk_\R(G)=\dim(M)-1$}\label{par:list2b}

In case (2-b) of the Main Theorem, the group $\Gamma$ is a lattice in a rank $n-1$ 
almost simple Lie group, and $\Gamma$ virtually embeds into $\Aut(M)^0$. This
implies that $\Aut(M)^0$ contains an almost simple complex Lie group $H$, the rank 
of which is equal to $n-1$. The goal in this section is to list all possible examples. 
Thus, our assumptions are 
\begin{enumerate}
\item[(i)] $H$ is an almost simple complex Lie group, and its rank is equal to $n-1$;
\item[(ii)] $M$ is a connected, compact, complex manifold and $\dim_\C(M)=n\geq 3$;
\item[(iii)] $H$ is contained in $\Aut(M)^0$.
\end{enumerate}
We now list all such possible pairs $(M,H)$. 

\begin{eg}\label{eg:homogeneous}
The group $\SL_n(\C)$ acts on $\P^{n-1}(\C)$ by linear projective transformations. 
In particular, $\SL_n(\C)$ acts on products of type $\P^{n-1}(\C)\times B$ where
$B$ is any Riemann surface. 

The action of  $\SL_n(\C)$ on $\P^{n-1}(\C)$ lifts to an action on the total space of the line bundles
${\mathcal{O}}(k)$ for every $k\geq 0$; sections of ${\mathcal{O}}(k)$ are in one-to-one
correspondence with homogeneous polynomials of degree $k,$ and the action of
$\SL_n(\C)$ on $H^0(\P^{n-1}(\C), {\mathcal{O}}(k))$ is the usual action on homogeneous 
polynomials in $n$ variables. 
Let $p$ be a positive integer and $E$ the vector bundle of rank $2$ over $\P^{n-1}(\C)$
defined by $E={\mathcal{O}}\oplus {\mathcal{O}}(p).$ Then $\SL_n(\C)$ acts
on $E,$ by isomorphisms of vector bundles. From this we get an action on the projectivized
bundle $\P(E),$ i.e. on a compact K\"ahler manifold $M$ which fibers over $\P^{n-1}(\C)$
with rational curves as fibers. 

When $k=1$, one can blow down the section of $\P(E)$ given by the line bundle ${\mathcal{O}}(1)$.
This provides a new smooth manifold with an action of $\SL_n(\C)$ (for other values of $k$, a
singularity appears). In that case, $\SL_n(\C)$ has an open orbit $\Omega$, the complement of
which is the union of a point and a smooth hypersurface $\P^{n-1}(\C)$.

A similar example is obtained from the $\C^*$-bundle associated to ${\mathcal{O}}(k).$
Let $\lambda$ be a complex number with modulus different from $0$ and $1.$ The quotient 
of this $\C^*$-bundle by multiplication by $\lambda$ along the fibers is a compact K\"ahler manifold,
with the structure of a torus principal bundle over $\P^{n-1}(\C).$ Since multiplication by 
$\lambda$ commutes with the $\SL_n(\C)$-action on  ${\mathcal{O}}(k),$ we obtain a (transitive) 
action of $\SL_n(\C)$ on this manifold. In this case, $M$ is not K\"ahler; if $k=1$, $M$ 
is the Hopf manifold, i.e. the quotient of  $\C^n\setminus\{0\}$ by the multiplication by $\lambda$. \end{eg}

\begin{eg}\label{eg:quadrics}
Let $H$ be the group $\SO_5(\C)$ (resp. $\SO_6(\C)$). The rank of $H$ is equal to $2$ (resp. $3$). 
The projective quadric $Q_3\subset \P^4(\C)$ (resp. $Q_4\subset \P^5(\C)$) 
given by the equation $\sum x_i^2=0$ is $H$-invariant, and has dimension $3$ (resp. $4$). 
The space of isotropic lines contained in $Q_3$ is parametrized by $\P^3(\C)$, so that $\P^3(\C)$
is a $\SO_5(\C)$-homogeneous space: This comes from the isogeny between $\SO_5(\C)$ and
$\Sp_4(\C)$ (see \cite{Fulton-Harris:book}, page 278), and provides another homogeneous 
space of dimension $\rk(\SO_5(\C))+1$. 

Similarly, $\SO_6(\C)$ is isogenous to $\SL_4(\C)$, and $\PSO_6(\C)$ acts transitively on $\P^3(\C)$. However,
 in this case, the rank of the group is equal to the dimension of the space (as in
Theorem \ref{thm:homogeneous}). 
\end{eg}

\begin{thm}\label{thm:list2b}
Let $M$ be a connected compact complex manifold of dimen\-sion~$n\geq 3$.
Let $H$ be an almost simple complex Lie group with $\rk_\C(H)=n-1$.  
If there exists an injective morphism $H\to \Aut(M)^0,$ then
$M$ is one of the following: 
\begin{itemize}
\item[(1)] a projective bundle $\P(E)$ for some rank $2$ vector bundle $E$ over $\P^{n-1}(\C),$ and then 
$H$ is isogenous to $\PGL_n(\C)$;

\item[(2)] a principal torus bundle over $\P^{n-1}(\C),$ and $H$ is isogenous to $\PGL_n(\C)$;

\item[(3)] a product of $\P^{n-1}(\C)$ with a curve $B$ of genus $g(B)\geq 2$, and then 
$H$ is isogenous to $\PGL_n(\C)$;

\item[(4)] the projective space $\P^n(\C),$ and $H$ is isogenous to
 $\PGL_n(\C)$ or to $\PSO_5(\C)$ when $n=3$;

\item[(5)] a smooth quadric of dimension $3$ or $4$ and $H$ is isogenous to 
$\SO_5(\C)$ or to $\SO_6(\C)$ respectively.
\end{itemize}
\end{thm}

The proof splits into three cases, according to the size of the orbits of 
$H$.

\subsubsection{Transitive actions}

Let us come back to the rank/dimension inequality obtained in Theorem \ref{thm:homogeneous}. 
Let $M$ be a connected compact complex manifold on  which a complex semi-simple Lie group
$S$ acts holomorphically and faithfully. Let $K\subset S$ be a maximal compact subgroup and let
$m$ be a point of $M$. Then
\begin{enumerate}
\item $\dim_\C(S)=\dim_\R(K)$;
\item $\dim_\R(K)=\dim_\R (K(m)) + \dim_\R(K_m)$ where $K(m)$ is the orbit of $m$ and $K_m$
is its stabilizer;
\item $K_m$ embeds into a maximal compact subgroup of $\GL(T_mM)$; in other words, $K_m$
is a closed subgroup of the unitary group $\U_n$, $n=\dim(M)$.
\end{enumerate}

The inequality 
\begin{equation}\label{eq:rkdim}
\dim_\C(S)\leq \dim_\R(M)+\dim_\R \U_n = 2n + n^2
\end{equation}
follows. Moreover, if the rank of $S$ is less than $n$, then $K_m$ has positive codimension in 
$\U_n$; this implies that ${\text{codim}}_\R(K_m)\geq 2n-2$ by classification of maximal subgroups
of $U_n$ or an argument similar to \S  \ref{part:homorank}. The inequality (\ref{eq:rkdim}) can therefore be strengthen, and
gives
\[
\dim_\C(S)\leq  n^2 + 2. 
\]
We now apply this inequality to the proof of Theorem \ref{thm:list2b} in case $H$ acts transtively. 
Thus, the semi-simple group $S$ is now replaced by the almost simple complex Lie group $H$, with  rank $r = n-1$. 

If the Lie algebra of $H$ is of type $B_r$ or $C_r$, i.e. $H$ is locally isomorphic to 
$\SO_{2r+1}(\C)$ or $\Sp_{2r}(\C)$, we have 
\[
\dim_\C(H)=  2r^2 + r \leq (r+1)^2 + 2  
\]
and thus $r^2 \leq r + 3$. This implies $r\leq 3$. 
When $r=2$, the group $H$ is locally isomorphic to $\SO_5(\C)$ and $\Sp_4(\C)$; there
are two examples of compact quotients of dimension $3$: The quadric $Q\subset \P^4(\C)$, and the projective
space $\P^3(\C)$ parametrizing the set of lines
contained in this quadric (see example \ref{eg:quadrics}).
When $r=3$, the group $H$ is isogenous to $\SO_7(\C)$ (resp. to $\Sp_6(\C)$) and there is no example of 
$H$-homogeneous compact complex manifold of dimension $4$ (see example \ref{eg:quadrics} and
\cite{Onishchik-Vinberg:EMS3}, page 169, 
\cite{Akhiezer:book}, page 65).

Let us now assume that $H$ is of type $D_r$, i.e. $H$ is isogenous to $\SO_{2r}(\C)$, with $r\geq 3$. 
We get  $ r^2 \leq 3r + 3$, so that $r=3$ and $H$ is isogenous to $\SO_6(\C)$. 
There is a unique homogeneous space $M$ of dimension $4$ for this group, namely the quadric
$Q\subset \P^5(\C)$. 

 Similarly, the inequality excludes the five exceptional groups ${\sf{E}}_6(\C),$ ${\sf{E}}_7(\C),$ ${\sf{E}}_8(\C),$
${\sf{F}}_4(\C)$, and ${\sf{G}}_2(\C)$: None of them acts transitively on a compact complex manifold of dimension $\rk(H)+1$.

The remaining case concerns the group $H=\SL_n(\C)$, acting transitively on a compact complex manifold 
$M$ of dimension $n\geq 3$. Write $M=H/L$ where $L$ is a closed subgroup of $H$. Two cases may occur: Either
$L$ is parabolic, or not. 

If $L$ is parabolic, then $M$ is a flag manifold of dimension $n$ for $\SL_n(\C)$. Flag manifolds for $\SL_n(\C)$ are well 
known, and only two examples satisfy our constraints. The first one 
is given by  the incidence variety $F\subset \P^2(\C)\times \P^2(\C)^\vee$  of pairs $(x,l)$ where $x$ is a point contained in the line $l$, 
or equivalently the set of complete flags of $\C^3$: This is a homogeneous space under the natural action 
of $\PGL_3(\C)$ and, at the same time, this is a $\P^1(\C)$-bundle over $\P^2(\C)$. 
The second example is given by the Grassmannian ${\mathbb{G}}(1,3)$ of lines
 in $\P^3$: This space has dimension $4$ and is homogeneous under the natural action of
$\PGL_4(\C)$. This example appears for the second time: By Pl\"ucker embedding, ${\mathbb{G}}(1,3)$ is a smooth quadric
in $\P^5(\C)$ and, as such, is a homogeneous space for $\SO_5(\C)$ (the groups $\SO_5(\C)$ and $\SL_4(\C)$ 
are isogenous, see page 286 of \cite{Fulton-Harris:book}).

If the group $L$ is not parabolic, then $L$ is contained in a parabolic group $P$ with $\dim(P)>\dim(L)$. 
This gives rise to a $H$-equivariant fibration 
\[
M\to H/P
\]
with $\dim(H/P)<n$. By Theorem \ref{thm:homogeneous}, $H/P$ is the projective space $\P^{n-1}(\C)$ and 
$\dim(P)=\dim(L)+1$. The fibers of the projection $M\to H/P$ are quotient of a one parameter group by a discrete
subgroup and, as such, are elliptic curves.  This implies that 
$M$ is an elliptic fibre bundle over $\P^{n-1}(\C)$, as in example \ref{eg:homogeneous}.

\subsubsection{Almost homogeneous examples}\label{part:ahe3}

Let us now assume that $M$ is not homogeneous under the action of $H$, but that $H$ has an open orbit 
$\Omega=H/L$; let $Z=M\setminus \Omega$ be its complement; this set is analytic and $H$-invariant. A theorem due to Borel  (\cite{Borel:1953})
asserts that the number of connected components of $Z$ is at most $2$.  By Proposition \ref{pro:ias}, each component of $Z$ is either a point
or a copy of $\P^{n-1}(\C)$; if one component is isomorphic to $\P^{n-1}(\C)$ then $H$ is isogenous to $\SL_n(\C)$ 
and acts transitively on this component.  Assume now that $Z$ contains an isolated point $m$. This point is fixed by the action 
of $H$, and this action can be linearized locally around $m$. Since $H$ has rank $n-1$ and $M$ has dimension $n$, the group $H$ is 
isogenous to $\SL_n(\C)$. Blowing up the point $m$, we replace $m$ by a copy of $\P^{n-1}(\C)$. 
Thus, $H$ is isogenous to $\SL_n(\C)$, and blowing up the isolated points of $Z$, we can assume that $Z$ is the union 
of one or two disjoint copies of $\P^{n-1}(\C)$ on which $H$ acts transitively. 
This situation has been studied in details in \cite{Huckleberry-Snow:1982} and \cite{Gilligan-Huckleberry:1981}; we now 
describe the conclusions of  \cite{Huckleberry-Snow:1982} and \cite{Gilligan-Huckleberry:1981} without proof.

Let $P$ be a maximal parabolic subgroup with $H/P=\P^{n-1}(\C)$ ($P$ is unique up to conjugacy). 

Suppose, first, that $Z$ is connected. Then $L\subset P$ (up to conjugacy), $M$ is a projective
rational manifold and it fibers equivariantly on $\P^{n-1}(\C)=H/P$; the fibers are isomorphic to $\P^1(\C)$, each of them intersecting $Z$
in one point (see \cite{Huckleberry-Snow:1982}). The intersection of each fiber with $\Omega$ is isomorphic to $\C$ and, at the same time, is
isomorphic to $P/L$; this is not possible for $n>2$ because all morphisms from the maximal parabolic group $P$ to the group 
$\Aff(\C)$ of holomorphic diffeomorphisms of $\C$ factor through the natural projection $P\to \C^*$, and there is no transitive action of
$\C^*$ on $\C$. 

Thus, $Z$ has indeed two connected components, as in  \cite{Gilligan-Huckleberry:1981} (see also \cite{Huckleberry-Snow:1982}). 
This case corresponds to $\P^1$-bundles over $\P^{n-1}(\C)$, 
as in example \ref{eg:homogeneous}:
$M$  fibers equivariantly on $\P^{n-1}(\C)$ with one dimensional fibers $F\simeq \P^1(\C)$, 
each of them intersecting $Z$ in two points; the two connected components of $Z$ are two sections
of the projection onto $\P^{n-1}(\C)$, which correspond to the two line bundles ${\mathcal{O}}$ and
${\mathcal{O}}(k)$ from example \ref{eg:homogeneous}. 

If $k=1$, one of the sections can be blown down to a fixed point (this process inverses the blow up construction described at the beginning 
of Section  \ref{part:ahe3}).

\subsubsection{No open orbit}

Let us now assume that $H$ does not have any open orbit. Then, blowing up all fixed points of $H$, 
all orbits have dimension $n-1$. By Theorem \ref{thm:homogeneous}, 
$H$ is isogenous to $\SL_n(\C)$ and its orbits are copies of $\P^{n-1}(\C)$.
In that case, the orbits define a locally trivial fibration of $M$ over a curve $B$. 
Let $A$ be the diagonal subgroup of $\SL_n(\C)$. The set of fixed points of $A$ 
defines $n$ sections of the fibration $M\to B$. This shows that this fibration is
trivial and $M$ is a product $\P^{n-1}(\C)\times B$. A posteriori, $H$ had no fixed point on $M$. 

%
%
\section{Invariant cones for lattices and Lie groups}\label{PART:4}
%
%

This paragraph contains preliminary results towards the proof of the Main Theorem in case (2-b).
Under the assumption of assertion (2-b), Proposition \ref{pro:ExtensionCohomology} applies, and
one can extend the action of $\Gamma$ on $W=H^{1,1}(M,\R)$ to an action of 
$G$; unfortunately, the nef cone $\Kabar(M)$ is not $G$-invariant a priori. In this section, 
we find a $G$-invariant subcone which is contained in $\Kabar(M)$. This is done 
in the general context of a linear representation of a semi-simple Lie group $G$, for which a lattice
$\Gamma\subset G$ preserves a salient cone.

\subsection{Proximal elements, proximal groups, and representations}

\subsubsection{Proximal elements and proximal groups} Let $V$ be a real vector space of finite dimension $k$. Let $g$ be an element of $\GL(V)$.
Let $\lambda_1(g) \geq \lambda_2(g) \geq ... \geq \lambda_k(g)$ be the moduli of the eigenvalues
of $g$, repeated according to their multiplicities. One says that $g$ is {\bf{proximal}} 
if $\lambda_1(g)>\lambda_2(g)$; in this
case, $g$ has a unique attracting fixed point $x^+_g$ in $\P(V)$. A subgroup  of $\GL(V)$ 
is proximal if it contains a proximal element, and a representation $G\to \GL(V)$ is proximal if its
image is a proximal subgroup.

If $\Gamma$ is a proximal subgroup of $\GL(V),$ the {\bf{limit set}} $\Lambda_\Gamma^{\P}$ of the
group $\Gamma$ in $\P(V)$ is defined as the closure of
the set $\{ x^+_g \vert g \in \Gamma, \, \, {\text{g is proximal}}\}$. 

\subsubsection{Proximal representations and highest weight vectors}
Let $G$ be a semi-simple Lie group and $A$ be a Cartan subgroup in $G$; 
let $\g$ and $\aa$ be their respective Lie algebras, and $\Sigma$ 
the system of restricted roots: By definition $\Sigma$ is the set of non-zero weights
for the adjoint action of $\aa$ on $\g$. One chooses a system of positive roots $\Sigma^+$. 
A scalar product $\langle \cdot \vert \cdot \rangle$
on $\aa$ is also chosen, in such a way that it is invariant by the Weyl group. 
One denotes by ${\sf{Wt}}$ the set of {\bf{weights}} of $\Sigma$; by definition 
\[
{\sf{Wt}}=\left\{   \lambda \in \aa \, \vert \, \,   \forall \, \alpha \in \Sigma, \, \,  \,  2\frac{\langle \lambda \vert \alpha \rangle}{\langle \alpha \vert \alpha \rangle}\in \Z \right\}.
\]
The set of  {\bf{dominant weights}} is ${\sf{Wt}}^+= \{\lambda \in {\sf{Wt}}\, \vert \, \, \forall \, \alpha \in \Sigma^+, \, \,  \, \langle \lambda\vert \alpha\rangle \geq 0 \}$. This set
corresponds to positive elements for the order defined on ${\sf{Wt}}$ by $\lambda\geq \lambda'$ if and only if  $\langle \lambda\vert \alpha\rangle \geq\langle \lambda'\vert \alpha\rangle$
for all $\alpha$ in $\Sigma^+$.

Let $\rho\colon G\to \GL(V)$ be an irreducible representation of $G$. This provides a representation of the Lie algebras $\g$
and $\aa$. By definition, the weights of $\aa$ in $V$ are the (restricted) weights of $\rho$.
This finite set has a maximal element
$\lambda$ for the order defined on ${\sf{Wt}}$: This {\bf{highest weight}} $\lambda$ is contained in ${\sf{Wt}}^+$, is unique, and determines 
the representation $\rho$ up to isomorphism.

The image of $G$ in $\GL(V)$ is proximal if and only if the eigenspace of $A$ corresponding to the highest weight $\lambda$
has dimension $1$ (see \cite{Abels-Margulis-Soifer:1995}).

If one starts with a representation $\rho$ which is not irreducible, one first splits it as a direct sum of irreducible
factors, and then apply the previous description to each of them; this gives a list of highest weights, one for each irreducible
factor. The maximal element in this list is the highest weight of $\rho$ (see \S~\ref{par:Reduc}).  

\subsection{Invariant cones}\label{par:InvCone}
In this paragraph we prove the following proposition.

\begin{pro}\label{pro:InvCone}
Let $\Gamma$ be a lattice in a connected semi-simple Lie group $G$.
Let $G\to \GL(V)$ be a real, finite dimensional, linear representation of
$G$. If $\Gamma$ preserves a salient cone $\Omega\subset V$ which is not reduced
to~$\{0\}$, the cone
$\Omega$ contains a $G$-invariant salient subcone which is not reduced to~$\{0\}$.
\end{pro}

Let $G$ be a connected semi-simple Lie group and $\Gamma$ be a Zariski dense subgroup of $G$.
Let $\rho:G\to \GL(V)$ be a real, finite dimensional, linear representation of $G$.  Assume
that $\rho(\Gamma)$ preserves a salient cone $\Omega$ with $\Omega \neq \{0\}$.
If the interior of $\Omega$ is empty, then $\Omega$ spans a proper $\Gamma$-invariant subspace of $V$;
since $\Gamma$ is
Zariski dense this proper invariant subspace is $G$-invariant. We can therefore restrict the study to this invariant 
subspace and assume that the interior of $\Omega$ is non empty. 

\begin{rem}
As the proof will show, if the action of $\Gamma$ on the linear span of $\Omega$ is not trivial, 
the action of $G$ on the linear span of its invariant subcone is also not trivial. In particular, if $G$
is simple, this action is faithful.
\end{rem}

\subsubsection{Irreducible representations}  

We first assume that $\rho$ is irreducible. 
Proposition 3.1, page 164, of \cite{Benoist:2000}, implies
that $\rho(\Gamma)$ is a proximal subgroup of $\GL(V)$, and the limit set $\Lambda_{\rho(\Gamma)}^\P$ of $\Gamma$ is contained in $\P(\overline{\Omega})$.
As a consequence, $\rho$ is a proximal representation of the group $G$ and the limit
set $\Lambda_{\rho(G)}^\P$ of $\rho(G)$ coincides with the orbit of the
highest weight line of its Cartan subgroup: This orbit is the unique closed orbit of $\rho(G)$ in $\P(V)$. 
As such, $\Lambda_{\rho(G)}^\P$ is a homogeneous space $G/P$, where $P$ is a parabolic subgroup of~$G$. 

Assume now that $\Gamma$ is a lattice in $G$. 
By \cite{Mostow:book}, lemma 8.5, all orbits $\Gamma\cdot x$ of $\Gamma$ in $G/P$ are dense, so that $\Lambda_{\rho(G)}^\P=G/P$ coincides
with $\Lambda_{\rho(\Gamma)}^\P$. In particular, $\Lambda_{\rho(G)}^\P$ is a $\rho(G)$-invariant 
subset of $\P(\overline{\Omega})$. The convex cone generated by  $\Lambda_{\rho(G)}^\P$ is
a closed and
$G$-invariant subcone of $\overline{\Omega}$. This proves Proposition \ref{pro:InvCone} for irreducible representations.

\subsubsection{General case}\label{par:Reduc}

Let us now consider a linear representation  $\rho:G\to \GL(V)$ which is not assumed
to be irreducible. Prasad and Raghunathan proved in \cite{Prasad-Raghunathan:1972} that $\Gamma$ intersects a 
conjugate of the Cartan subgroup $A'\subset G$ on a cocompact lattice $A_\Gamma' \subset A'$. Changing
$A$ into $A'$, we assume that $\Gamma$ intersects $A$ on such a lattice~$A_\Gamma$.

Since $G$ is semi-simple, $V$ splits into a direct sum of irreducible factors; let $\lambda$ be
the highest weight of $(\rho,V)$, let $V_1$, ... $V_m$ be the irreducible factors corresponding to this weight, 
and let $V'$ be the direct sum of the $V_i$:
\[
V':=\bigoplus_{1\leq i \leq m} V_i.
\] 
By construction, all representations $V_i$, ${1\leq i \leq m}$, are isomorphic.
\begin{lem}\label{lem:4.2}
Since $\Gamma$ is a lattice, $\overline{\Omega}$ intersects the sum $V'$
of the highest weight factors on a closed, salient cone $\Omega'$ which is not reduced to zero. 
\end{lem}

\begin{proof}
If $u$ is any element of $\Omega$, one can decompose
$u$ as a sum $\sum_\chi u_\chi$ where each $u_\chi$ is an eigenvector of the Cartan subgroup $A$ corresponding
to the weight $\chi$. Since $\Omega$ has non empty interior, we can choose such an element $u$ with a 
non zero component $u_\lambda$ for the highest weight~$\lambda$. Since $A_\Gamma$ is a lattice in $A$,
there is a sequence of elements $\gamma_n$ in $A_\Gamma$ such that 
\[
\frac{\gamma_n(u)}{\Vert \gamma_n(u)\Vert}=\sum_\chi \frac{\chi(\gamma_n)}{\Vert \gamma_n(u)\Vert} u_\chi
\]
converges to a non zero multiple of $u_\lambda$. Since $\Omega$ is $\Gamma$-invariant and all $\gamma_n$
are in $\Gamma$, we deduce that $\overline\Omega$ intersects $V'$.
\end{proof}

The subspace of $V'$ which is spanned by $ \Omega'$
is a direct sum of highest weight factors; for simplicity, we can therefore assume that $V'$ is spanned by $\Omega'$.
In particular, the interior $\Int(\Omega')$ is a non-empty subset of $V'$.

Let $\pi_i$ be the projection of $V'=\bigoplus V_i$ onto the factor $V_i$. The image of $\Int(\Omega')$ by $\pi_1$ is
an open subcone   $\pi_1(\Omega')$ in $V_1$. 

 If this cone is salient, the previous paragraph shows that the
representation $(\rho_1, V_1)$ is proximal. Thus, all $V_i$ can be identified to a unique proximal representation $R$, with a 
given highest weight line $L=\R u^+$. We obtain $m$ copies $L_i$ of $L$, one in each copy $V_i$ of $R$. Apply Lemma 
\ref{lem:4.2} and its proof: Since $\Omega'$ 
is $\Gamma$-invariant, $\Gamma$ is a lattice, and $\Omega'$ has non empty interior, there is a point $v\in L_1\oplus ... \oplus L_m$
which is contained in $\Omega'$. Let $(a_1, ..., a_m)$ be the real numbers such that $v=(a_1u^+, ..., a_mu^+)$. The
diagonal embedding $R\to V'$, $w\mapsto (a_1w, ..., a_m w)$ determines an irreducible sub-representation of $G$
into $V$ that intersects $\Omega'$, and the previous paragraph shows that $G$ preserves a salient
subcone of $\Omega'$.

If the cone $\pi_1(\Omega')$ is not salient, the fiber $\pi_1^{-1}(0)$ intersects $\Omega'$ 
on a $\Gamma$-invariant salient subcone; this reduce the number of irreducible factors from $m$ to $m-1$,
and enables us to prove Proposition  \ref{pro:InvCone} by induction on the number $m$ of factors $V_i$.

%
%

\section{Linear representations, ample classes and tori}\label{PART:5}
%
%

We now prove the Main Theorem. Recall that $G$ is a connected, almost simple, real Lie group with real rank $\rk_\R(G)\geq 2$,
that $A$ is a Cartan subgroup of $G$, and that $\Gamma$ is a lattice in $G$ acting on a connected compact K\"ahler manifold
$M$ of dimension $n$. 

From Section \ref{par:firstpart}, we know that the rank of $G$ is at most $n$ and, in case $\rk_\R(G)=n$, the group $G$ is isogenous
to $\SL_{n+1}(\R)$ or $\SL_{n+1}(\C)$ and $M$ is isomorphic to $\P^n(\C)$. We now assume that the rank of $G$ satisfies the next critical 
equality $\rk_\R(G)=n-1$. 
According to Proposition \ref{pro:ExtensionCohomology}, two possibilities can occur.
\begin{itemize}
\item The image of $\Gamma$ is virtually contained in $\Aut(M)^0$;  Theorem \ref{thm:list2b} in Section \ref{par:list2b} gives
the list of possible pairs $(M,G)$. This corresponds to assertion (2-a) in the Main Theorem.
\item The action of $\Gamma$ on the cohomology of $M$ is almost faithful and  virtually extends
to a linear representation of $G$ on $H^*(M,\R)$. 
\end{itemize}
Thus, in order to prove the Main Theorem, we replace $\Gamma$ by a finite index subgroup and assume that the action of $\Gamma$ on the
cohomology of $M$ is  faithful and extends to a linear representation of $G$.  Our aim is to prove that all such examples
are Kummer examples (when $\rk_\R(G)=\dim_\C(M)-1$).

We  denote by $W$ the space $H^{1,1}(M,\R)$, by $\lambda_W$ the highest weight of the representation 
$G\to \GL(W)$ and by $E$ the direct sum of the irreducible factors $V_i$ of $W$ corresponding to the 
weight $\lambda_W$ (all $V_i$ are isomorphic representations).

\subsection{Invariant cones in $\Kabar(M)$}

Since the K\"ahler cone $\Ka(M)$ is a $\Gamma$-invariant, convex, and salient cone in $W$
with non empty interior,
Proposition \ref{pro:InvCone} asserts that  $\Kabar(M)$ contains a non-trivial $G$-invariant subcone. More
precisely,  $\Kabar(M)\cap E$ contains a $G$-invariant salient subcone  $\Kabar_E$ which
is not reduced to~$\{0\}$, and the action of $G$ on the linear span of $\Kabar_E$ is faithful (see \S \ref{par:InvCone}).

From now on, we replace $E$ by the linear span of the cone $\Kabar_E$. Doing this, 
the cone $\Kabar_E$ has non empty interior in $E$, and is a $G$-invariant subcone of $\Kabar(M)$.
Since $G$ is almost simple, the representation $G\to \GL(E)$ is unimodular. 
Thus, the action of the Cartan subgroup $A$ on $E$ is unimodular, faithful and diagonalizable. 

\subsection{Actions of abelian groups}\label{par:AAG}

We now focus on a slightly more general situation, and use ideas from \cite{Dinh-Sibony:Duke}.
Let $A$ be the additive abelian group $\R^m$, with $m\geq 1$; in the following paragraph, 
$A$ will be a Cartan subgroup of $G$, and thus $m=\rk_\R(G)$ will be equal to $\dim(M)-1$. Let $E$ be a subspace of $W$ 
and $\Kabar_E$ be a subcone of $\Kabar(M)\cap E$ with non empty interior. Let $\rho$ be a continous
representation of $A$ into $\GL(H^*(M,\R))$ by cohomological automorphisms. Assume that 
\begin{itemize}
\item[(i)] $\rho(A)$ preserves $E$ and $\Kabar_E$;
\item[(ii)] the restriction $\rho_E\colon A\to \GL(E)$ is diagonalizable, unimodular, and faithful. 
\end{itemize}
From (ii), there is a basis of $E$ and morphisms $\lambda_i:A\to \R$, $1\leq i\leq \dim(E),$ 
such that the matrix of $\rho_E(a)$ in this basis is  diagonal, with diagonal coefficients $\exp(\lambda_i(a))$. 
The morphisms $\lambda_i$ are the weights of $\rho_E$; the set of weights 
\[
\Lambda= \left\{   \lambda_i,  1\leq i\leq \dim(E) \right\}
\]
is a finite subset of $A^\vee$ where $A^\vee$, the dual of $A$, is identified with the
space of linear forms on the real vector space $A=\R^m.$ 
The convex hull  of $\Lambda$  is a polytope ${\mathcal{C}}(\Lambda)\subset A^\vee$ and the set of its extremal vertices is a subset 
$\Lambda^+$ of $\Lambda$; equivalently, a weight $\lambda$ is {\bf{extremal}} if and only if there is an element
$a\in A$ such that 
\[
\lambda(a)>\alpha(a), \quad \forall \alpha\in \Lambda\setminus\{\lambda\}.
\]
Since any convex set is the convex hull of its extremal points, $\Lambda^+$ is not empty and ${\mathcal{C}}(\Lambda^+)$
coincides with ${\mathcal{C}}(\Lambda)$.

For all weights $\alpha\in \Lambda$, we denote by $E_\alpha$ the eigenspace of $A$ 
of weight $\alpha$:
\[
E_\alpha=\left\{ u\in E\, \vert  \quad  \forall a\in A, \, \, \rho_E(a)(u)=e^{\alpha(a)}u   \right\}.
\]
We denote by $E^+$ the vector subspace of $E$ which is spanned by the $E_\lambda$ where $\lambda$ describes
$\Lambda^+$.

\begin{lem}\label{lem:NTI} The following three properties are satisfied.
\begin{enumerate}
\item The representation $\rho_{E^+}\colon A\to \GL(E^+)$ is injective.
\item The convex hull ${\mathcal{C}}(\Lambda)$ of $\Lambda^+$ contains the origin in its interior; in particular 
the cardinal of $\Lambda^+$ satisfies $\vert \Lambda^+\vert \geq \dim(A)+1$.
\item For all $\lambda\in \Lambda^+$ we have $E_\lambda\cap \Kabar_E\neq \{0\}$.
\end{enumerate}
\end{lem}

\begin{proof} The three properties are well known. 

\vv

\noindent{Property (1).---} The kernel of $\rho_{E^+}$ is defined by the set of linear
equations $\lambda(a)=0$ where $\lambda$ describes $\Lambda^+$. Since all weights
$\alpha\in \Lambda$ are barycentric combinations of the extremal weights, the kernel 
of $\rho_{E^+}$ is contained in the kernel of $\rho_E$. Property (1) follows from the 
injectivity of $\rho_E$.

\vv

\noindent{Property (2).---} 
 If ${\mathcal{C}}(\Lambda)$ has empty interior, it is contained in a strict subspace
of $A^\vee$, contradicting Property (1). In particular, the cardinal of $\Lambda^+$ satifies $\vert \Lambda^+\vert \geq \dim(A)+1$.
Since the sum of all weights $\lambda_i(a)$, repeated with multiplicities, is the logarithm of the determinant of 
$\rho_E(a)$ and the representation is unimodular, this sum is $0$; in other words, the origin of $A^\vee$ is a barycentric
combination of all extremal weights with strictly positive masses. This shows that the origin is in the interior of ${\mathcal{C}}(\Lambda)$.

\vv

\noindent{Property (3).---} The proof is similar to the proof of Lemma \ref{lem:4.2}. Let $\lambda$ be an extremal weight and let $a\in A$ satisfy
$\lambda(a)>\alpha(a)$ for all $\alpha\in \Lambda\setminus \{\lambda\}$. Let $u$ be any element of $\Kabar_E$; write
$u$ as a linear combination 
$
u=\sum_{\alpha\in \Lambda} u_\alpha
$
where $u_\alpha\in E_\alpha$ for all $\alpha$ in $\Lambda$. Since $\Kabar_E$ has non empty
interior, we can choose $u$ in such a way that $u_\lambda\neq 0$. Then the sequence
\[
\frac{\rho_E(na)(u)}{\exp(n\lambda(a))} 
\]
is a sequence of elements of $\Kabar_E$ that converges towards $u_\lambda$ when $n$ goes to $+\infty$. Since $\Kabar_E$
is closed, property (3) is proved.
\end{proof}

\begin{lem}\label{lem:NTP}
Let $k$ be an integer satisfying 
\[
1 \leq k \leq \min(\dim(M), \dim(A)+1).
\]
Let $\lambda_i\in \Lambda$, $1\leq i\leq \dim(A)+1$, be  distinct weights, and
$w_i$ be non zero elements in $E_{\lambda_i}\cap\Kabar_E$. 
For all multi-indices $I=(i_1, ... i_k)$ of distinct integers $i_j\in \{1, \ldots, \dim(A)+1\}$ the wedge product
\[
w_I:=w_{i_1}\wedge \ldots \wedge w_{i_k}
\] 
is different from $0$.
\end{lem}

The proof makes use of the following proposition which is due to Dinh and Sibony (see \cite{Dinh-Sibony:Duke}, Corollary 3.3 and Lemma 4.4).
Lemma 4.4 of \cite{Dinh-Sibony:Duke} is stated for cohomological automorphisms that are
induced by automorphisms of $M$, but the proof given in \cite{Dinh-Sibony:Duke} extends to all cohomological automorphisms.

\begin{pro}\label{pro:DS} Let $M$ be a connected compact K\"ahler manifold.
Let $u$ and $v$ be elements of $\Kabar(M)$. 
\begin{enumerate}
\item If $u$ and $v$ are not colinear, then $u\wedge v \neq 0$.
\item Let $v_1, \ldots, v_l$, $l\leq n-2$, be elements of $\Kabar(M)$. If 
$v_1\wedge \ldots \wedge v_l\wedge u$ and $ v_1\wedge \ldots \wedge v_l\wedge v$
are non zero eigenvectors with distinct eigenvalues for a cohomological automorphism, then 
$(v_1\wedge \ldots \wedge v_l)\wedge (u\wedge v) \neq 0.
$
\end{enumerate}

\end{pro}

\begin{proof}[Proof of Lemma \ref{lem:NTP}]
The proof is an induction on $k$. Since all $w_i$ are assumed to be different from $0$, 
the property is established for $k=1$. Assume that the property holds for all multi-indices $I$ of
length $\vert I\vert =k$ with 
\[
k\leq \min(\dim(M),\dim(A)+1)-1.
\] 
Let $J=(i_1, \ldots, i_{k+1})$ be a multi-index of 
length $k+1$. Let $v_1, \ldots, v_{k-1}$ denote the vectors $w_{i_1}, \ldots, w_{i_{k-1}}$, let $u$ 
be equal to $w_{i_k}$ and $v$ be equal to $w_{i_{k+1}}$. Since the property is proved for length $k$,
we know that 
\[
v_1\wedge \ldots \wedge v_l\wedge u \quad {\text{and}} \quad v_1\wedge \ldots \wedge v_l\wedge v
\] 
are two non zero eigenvectors of $A$ with respective weights 
\[
\lambda_{i_{k}}+ \sum_{j=1}^{j=k-1} \lambda_{i_j} \quad {\text{and}} \quad
\lambda_{i_{k+1}}+ \sum_{j=1}^{j=k-1} \lambda_{i_j}.
\]
These two weights are different because $\lambda_{i_{k}}\neq \lambda_{i_{k+1}}$.
Thus, property (2) of proposition \ref{pro:DS}  can be applied, and it implies that $w_J$ 
is different from zero. The Lemma follows by induction.
\end{proof}

Let us now assume that $\dim_\R(A)=\dim_\C(M)-1$, i.e. $m=n-1$.  
According to property (2) in Lemma \ref{lem:NTI}, we can find $n=\dim(A)+1$ 
extremal weights $\lambda_i$ such that all linear maps
\[
a\mapsto (\lambda_1(a), \ldots , \lambda_{i-1}(a), \lambda_{i+1}(a), \ldots,  \lambda_{n}(a)), \quad 1\leq i \leq n,
\]
are bijections from $A$ to $\R^{n-1}$.
B property (3) in Lemma \ref{lem:NTI}, there exist elements $w_i$ in $E_{\lambda_i}\cap\Kabar_E\setminus\{ 0\}$ for all $1\leq i \leq \dim(A)+1$.
Once such a choice of vectors $w_i$ is made, we define $w_A$ as the sum 
\[
w_A=w_1+w_2+\ldots + w_{n}.
\]
This class is nef and, by Lemma \ref{lem:NTP}, its $n$-th power is different from zero ; it is a sum 
of products $w_{i_1}\wedge ... \wedge w_{i_n}$ which are positive because all classes $w_i$ are nef. 
Thus,  
\[
w_A^{\wedge n}= w_A\wedge w_A \wedge \ldots \wedge w_A >0.
\]
According to Section \ref{par:NefBig}, this proves the following corollary.

\begin{cor}\label{cor:NB}
 If $\dim_\R(A)=\dim_\C(M)-1$, the class $w_A$ is nef and big.
\end{cor}




\subsection{A characterization of torus examples}\label{par:toresc1}

Let us apply the previous paragraph to the Cartan subgroup $A$ of $G$; by assumption, 
$G$ has rank $n-1$ and thus $\dim_\R(A)=\dim_\C(M)-1=n-1$. The groups $G$ and
$A$ act on $W$ and preserve $E$, and we denote by $\rho_E(g)$ the endomorphism 
of $E$ obtained by the action of $g\in G$ on $E$ (thus, $\rho_E(g)$ is the restriction of $g^*$
if $g$ is in $\Gamma$).
We keep the notation of Section \ref{par:AAG}, as well as the choice of classes $w_i$ and $w_A$.
According to Corollary \ref{cor:NB}, the class $w_A$ is nef and big. 

\begin{pro}\label{pro:toresc1}
If the class $w_A$ is a K\"ahler class then, up to a finite cover, $M$ is a torus.
\end{pro}

\begin{rem}
In Section \ref{par:OBAIS}, this result is applied in the slightly more general context where $M$ is an
orbifold with isolated singularities.
\end{rem}

\begin{proof}
Let $c_1(M)\in H^{1,1}(M,\R)$ and $c_2(M)\in H^{2,2}(M,\R)$ be the first and second Chern classes of $M$. 
Both of them are invariant under the action of $\Gamma$, and therefore also under the action of $G$. 

Let $u\in W$ be a $G$-invariant cohomology class. Let $I=(i_1, ..., i_{n-1})$ be a multi-index
of length $n-1$, and $w_I$ be the product $w_{i_1}\wedge \ldots \wedge w_{i_{n-1}}$. Let $v$ 
be the class of type $(n,n)$ defined by $v=w_I\wedge u$. Since $u$ is $A$-invariant we have 
\[
\rho_E(a)(v)=\exp\left( \sum_{j=1}^{n-1} \lambda_{i_j}(a)\right) v.
\]
Since $v$ is an element of $H^{n,n}(M,\R)$ and the action of $G$ is trivial on
$H^{n,n}(M,\R)$, we get the alternative: Either $v=0$ or  $\sum_{j=1}^{j=n-1} \lambda_{i_j}(a)=0$ for all $a\in A$. 
Thus,  property (i) of the extremal weights $\lambda_i$ implies that $v$ is equal to $0$
for all choices of multi-indices $I$ of length $n-1$. 
As a consequence, $u$ is a primitive element with respect to the K\"ahler class $w_A$:
\[
\int_M w_A^{n-1}\wedge u = 0.
\] 
In the same way, one proves that $w_A^{n-2}\wedge u=0$ for all $G$-invariant cohomology classes  
$u$ in $H^{2,2}(M,\R)$. 

Let us apply this remark to the first Chern class $c_1(M)$. Since this class is invariant, it is
primitive with respect to $w_A$. Since $c_1(M)^2$ is also $G$-invariant, 
\[
w_A^{n-2}\wedge c_1(M)^2 = 0\, ;
\]
From Hodge index theorem we deduce that  $c_1(M)=0$. Yau's theorem provides a Ricci flat
K\"ahler  metric on $M$ with K\"ahler form $w_A$, and Yau's formula reads
\[
\int_M w_A^{n-2}\wedge c_2(M)= \kappa \int_M \Vert {\sf{Rm}} \Vert^2 w_A^n
\] 
where ${\sf{Rm}}$ is the Riemannian tensor and $\kappa$ is a positive constant (see \cite{Besse:book}, page 80, and \cite{Kobayashi:book}, \S IV.4 page 112--118). From the invariance
of $c_2(M)$ we get $w_A^{n-2}\wedge c_2(M)=0$ and then ${\sf{Rm}}=0$. This means that $M$ is flat
and thus $M$ is finitely covered by a torus $A$.
\end{proof}

Using this proposition, we now need to  change the big and nef class $w_A$ into an ample class by a modification of $M$. 
This is the main goal of the following paragraph.

\subsection{Obstruction to ampleness, invariant subsets, and Kummer examples}\label{par:OBAIS}

\subsubsection{} Let us start with the following simple fact. 

\begin{pro}\label{pro:ruseanalytique}
Let $B$ be an irreducible real analytic subset of the vector space $H^{1,1}(M,\R)$. Assume that
\begin{itemize}
\item[(i)] all classes $w$ in $B$ are big and nef classes but
\item[(ii)] none of them is ample. 
\end{itemize}
Then there exists an integer $d$ with $0<d< n$ and a complex analytic subset 
$Y_0\subset M$ of dimension $d$ such that 
$
\int_{Y_0} w^d= 0
$
for all classes $w$ in $B$.
\end{pro}

\begin{proof}
The set of classes $[Y]$ of irreducible analytic subsets $Y\subset X$ is countable. For all
such classes $[Y]$, let $Z_{[Y]}$ be the closed, analytic subset of $B$ which 
is defined by
\[
Z_{[Y]}=\left\{ w \in B \, \vert \, \, \int_{Y}w^{\dim(Y)} = 0\right\}.
\]
Apply Section \ref{par:NefBig}. Since all elements of $B$ are nef and big but none of them is ample, 
the family of closed subsets $Z_{[Y]}$ with $\dim(Y)\geq 1$ covers $B$. By Baire's theorem, one of the subsets
$Z_{[Y]}$ has non empty interior. Let $Z_{[Y_0]}$ be such a subset, with $\dim(Y_0)\geq 1$. The map 
\[
w \mapsto \int_{Y_0} w^{\dim(Y_0)}
\]
is algebraic and vanishes identically on an open subset of $B$. Since $B$ is an irreducible
analytic subset of $H^{1,1}(M,\R)$, this map vanishes identically. 
\end{proof}

 \subsubsection{} Coming back to the proof of Theorem A, and assuming that $w_A$ is not ample 
 (for all Cartan subgroups $A$ of $G$), we consider the orbit of the class $w_A$ 
under the action of $G$. This orbit $B=G.w_A$ satisfies the following properties:
\begin{enumerate}
\item  $B$ is made of big and nef classes, but none of them is ample;
\item $B$ is a connected Zariski open subset in an irreducible algebraic subset of $E$.
\end{enumerate} 
We can thus apply Proposition \ref{pro:ruseanalytique} to the set $B$. Let $Z$ be the union 
of analytic subsets $Y\subset M$ such that $0<\dim(Y)<\dim(M)$ and 
\[
\int_Y w^{\dim(Y)}=0, \quad \forall w \in B.
\]
Proposition \ref{pro:ruseanalytique} and Section \ref{par:NefBig} show that $Z$ is a non empty 
proper analytic subset of $M$. Since $B$ is the orbit of $w_A$ under the action of $G$, this set
is $\Gamma$-invariant. 

\subsubsection{} Let  us now apply Theorem \ref{thm:ih} to the subset $Z\subset M$ that is constructed
in the previous paragraph. We  get a birational morphism $\pi:M\to M_0$ and conclude that 
the image of $w_A$ in $M_0$ is ample. From Section \ref{par:toresc1} and Proposition \ref{pro:toresc1} applied in the orbifold
context, we deduce that $M_0$ is covered 
by a torus $A$. 

Let us be more precise. In our case, $M_0$ is a connected
orbifold with trivial Chern classes $c_1(M_0)$ and $c_2(M_0).$ This implies that there 
is a flat K\"ahler metric on $M_0$ (see \cite{Kobayashi:book}).
The universal cover of $M_0$ (in the orbifold sense) is then isomorphic to 
$\C^n$  and the (orbifold) fundamental group $\pi_1^{orb}(M_0)$ acts 
by affine isometries on $\C^n$ for the standard euclidean metric. 
In other words, $\pi_1^{orb}(M_0)$
 is identified to a cristallographic group $\Delta$ of affine motions of $\C^n.$ 
 Let $\Delta^*$ be the group of translations contained in $\Delta.$ 
 Bieberbach's theorem shows that  (see \cite{Wolf:book}, chapter
 3, theorem 3.2.9).  
 \begin{itemize}
\item[a.-] $\Delta^*$ is a lattice in $\C^n$;
\item[b.-] $\Delta^*$ is the unique maximal and normal free abelian subgroup 
of $\Delta$ of rank $2n.$
 \end{itemize}
The torus $A$   is the quotient  of $\C^n$ by this group of translations. By construction, $A$ covers $M_0$.
 Let $F$ be the quotient  group $\Delta/\Delta^*$; we identify it to the 
 group of deck transformations of the covering $\epsilon:A\to M_0.$ To conlude the proof 
 of the Main Theorem, all we need to do is to lift virtually the action of $\Gamma$ on $M_0$ to an action
 on $A$. This is done in the following lemma.

\begin{lem}\label{lem:Fcyclic} $\,$
\begin{enumerate}
\item A finite index subgroup of $\Gamma$ lifts 
to $\Aut(A).$
\item Either $M_0$ is singular, or $M_0$ is a torus. 
\item If $M_0$ is singular, then $M_0$ is a quotient of the torus 
$A$ by a homothety $(x,y,z)\mapsto (\eta x, \eta y, \eta z),$ where $\eta$
is a root of $1.$
\end{enumerate}
\end{lem}

\begin{proof}
By property (b.) all automorphisms of $M_0$ lift to $A.$ 
Let $\overline{\Gamma} \subset \Aut(A)$ be the group of automorphisms of $A$ made of all possible lifts of elements
of~$\Gamma.$ So, $\overline{\Gamma}$ is an extension of  $\Gamma$ by the group $F$: 
\[
1\to F\to \overline{\Gamma} \to \Gamma \to 1.
\]
Let $L:\Aut(A)\to \GL_n(\C)$ be the morphism which applies each automorphism
$f$ of $A$ to its linear part $L(f).$ Since $A$ is obtained as the quotient of $\C^n$ 
by all translations contained in $\Delta,$ the restriction of $L$ to $F$ is injective. 
Let $N\subset \GL_n(\C)$ be the normalizer of $L(F)$.
The group $L( \overline{\Gamma})$ normalizes $L(F).$ Hence we have a well defined
morphism $\overline{\Gamma}\to N,$ and an induced morphism $\delta:\Gamma\to N/L(F)$. 
Changing $\Gamma$ into a finite index subgroup,  
$\delta$ is injective. 
Since $\Gamma$ is a lattice in an almost simple Lie group of rank $n-1$, 
the Lie algebra of $N/L(F)$ contains a subalgebra of rank $n-1$. Since  $\sll_n(\C)$ is the unique complex subalgebra of rank $n-1$ in $\gll_n(\C)$,
we conclude that $N$ contains $\SL_n(\C)$. It follows that $L(F)$ is contained 
in the center $\C^* {\text{Id}}$ of $\GL_n(\C)$. 

Either $F$ is trivial, and then $M_0$ coincides with the torus $A,$ or $F$ is
a cyclic subgroup of $\C^* {\text{Id}}.$ In the first case, there is no need to 
lift $\Gamma$ to $\Aut(A).$
In the second case, we fix a generator $g$ 
of $F,$ and denote by $\eta$ the root of unity such that $L(g)$ is the multiplication
by $\eta.$ The automorphism $g$ has at least one (isolated) fixed point $x_0$ 
in $A$. Changing $\Gamma$ into a finite
index subgroup $\Gamma_1,$ we can assume that $\overline{\Gamma_1}$ 
fixes $x_0.$ The linear part $L$ embeds $\overline{\Gamma_1}$ 
into $\GL_n(\C)$. Selberg's lemma assures that a finite index subgroup
of $\overline{\Gamma_1}$ has no torsion. This subgroup does not intersect
$F,$ hence projects bijectively onto a finite index subgroup of $\Gamma_1.$ 
This proves that a finite index subgroup $\Gamma_1$ of $\Gamma$ lifts to $\Aut(A).$
\end{proof}

%
%

\section{Classification of Kummer examples}\label{par:Kummer}
%
%

In this section, we list all Kummer examples of dimension $n\geq 3$ with an action of a lattice $\Gamma$
in a rank $n-1$ simple Lie group $G$, up to commensurability  and isogenies. 
The main step is to classify tori $\C^n/\Lambda$  such that $\Aut(\C^n/\Lambda)$ contains a copy of $\Gamma$. 
As seen in the proof of Lemma \ref{lem:Fcyclic}, a finite index subgroup of $\Gamma$ lifts to a linear representation 
into $\SL_n(\C)$ that preserves the lattice $\Lambda$. Margulis theorem implies that this linear representation 
virtually extends to a representation of $G$ itself. Thus, we have to list triples $(G,\Gamma, \Lambda)$ 
where $G$ is a real almost simple Lie group  represented in $\SL_n(\C)$, $\Gamma$ is a lattice in $G$, 
$\Lambda$ is a lattice in $\C^n$, and $\Gamma$ preserves $\Lambda$. This is done in paragraphs 
\ref{part6:prelim} to \ref{par:quater}: The list is up to commensurability for $\Gamma$, and up to
isogeny for $\C^n/\Lambda$. 
Then we discuss Kummer examples in paragraph \ref{part6:Kummer}.

\subsection{Preliminaries}\label{part6:prelim}

If a rank $n-1$ connected simple real Lie group $G$ acts on $\C^n$ linearly, then $G$ is locally isomorphic
to $\SL_n(\R)$ or $\SL_n(\C)$. We can therefore assume that $G$ is either 
$\SL_n(\C)$ or $\SL_n(\R)$ and  $\Gamma$ is a lattice in $G$.

For actions of lattices $\Gamma\subset \SL_n(\C)$ on tori,  proposition 8.2 of \cite{Cantat-Zeghib} can be applied: 
There is a negative integer $d$, such that $\Gamma$ is commensurable to $\SL_n({\mathcal{O}}_d)$ where 
${\mathcal{O}}_d$ is the ring of integers of the quadratic number field $\Q(\sqrt{d})$, and the torus $M$ is
isogeneous to $(\C/{\mathcal{O}}_d)^n$. 

We can therefore restrict our study to the case of $\SL_n(\R)$.

\subsubsection{Setting} 

In what follows, $\Gamma$ is a lattice in $G=\SL_n(\R)$, $G$ acts linearly on $V=\C^n$,
 by a morphism $\rho:G\to \SL_n(\C)$, and $\Gamma$ preserves a  lattice $\Lambda\subset V$.

If we forget the complex structure, we can identify the vector space $V=\C^n$ with $\R^{2n}$ and the lattice $\Lambda$ 
with $\Lambda=\Z^{2n}$; the complex structure on $V$ is then given by a linear operator $J\in \GL_{2n}(\R)$ with $J^2=-{\text{Id}}_{2n}$. The
linear representation $\rho$ of $G$ preserves the complex structure. As a consequence,
$\rho$  is equivalent to the diagonal representation on $P\times P$ where $P=\R^{n}$ is the standard representation of $G$. 
More precisely, the vector space $V$ splits as $P_1\oplus J(P_1)$ where $P_1$ is a $G$-invariant totally real $n$ dimensional
subspace; restricted to $P_1$ and to $J(P_1)$, the representation of $G$ is conjugate to its standard representation $P$.
The complex structure $J$ acts as follows: If $(u,v)$ is a point in $P\times P\simeq V$, then $J(u,v)=(-v,u)$.

Up to finite index, the lattice $\Gamma$ coincides with the lattice
\[
\{g \in G \, \vert \, \rho(g)(\Lambda)=\Lambda\},
\]
that is with the (preimage of the) intersection $\rho(G)\cap \SL_{2n}(\Z)$.
In particular, $\rho(G)\cap \SL_{2n}(\Z)$ is Zariski dense in $\rho(G)$. 

\subsubsection{Centralizer}
Let $C_G\subset \SL(V)$ be the centralizer of $\rho(G)$:
\[
C_G=\left\{ h \in \SL(V)\, \vert \, h\rho(g)=\rho(g)h\quad \forall \, g \in G     \right\}.
\]
 As a Lie group, $C_G$ is isomorphic to $\SL_2(\R)$ acting on $V=P\times P$ by 
\[
(u,v) \in P\times P \mapsto (au+bv, cu+dv).
\]
In particular, it does not preserve the complex structure $J$. 
Since $\rho(G)\cap \SL_{2n}(\Z)$ is Zariski dense in $\rho(G)$, the centralizer $C_G$ is defined
over $\Z$ in $\SL_{2n}(\R)$; hence $C_G(\Z)=C_G\cap \SL_{2n}(\Z)$ is an arithmetic
lattice in $C_G$. As such, either $C_G(\Z)$ is not cocompact and is then commensurable  
to $\SL_2(\Z)$, or $C_G(\Z)$ is cocompact and is then commensurable to a lattice
derived from a quaternion algebra (see \S \ref{par:quaternions} below, \cite{Katok:book}, and \cite{Witte:prebook}, chapter 6).

\begin{lem}\label{lem:pepr}
The following properties are equivalent. 
\begin{enumerate}
\item $\rho(G)$ preserves an $n$- dimensional plane $P'$ in $V$ which is defined over $\Q$.
\item $\rho(G)$ is conjugate by an element of $\SL_{2n}(\Q)$ to the standard diagonal
group $\{ (A,A)\, \vert \, A \in \SL_n(\R)\}$ in $\SL_{2n}(\R)$.
\item up to finite indices, $\rho(\Gamma)$ is conjugate to 
the diagonal copy of $\SL_n(\Z)$ in $\SL_{2n}(\R)$.
\item up to finite indices, $\Gamma$ is conjugate to $\SL_n(\Z)$ in $G=\SL_n(\R)$.
\item $C_G(\Z)$ is not cocompact, and is thus commensurable to $\SL_2(\Z)$.
\end{enumerate}
\end{lem}

\begin{proof}
If $\rho(G)$ preserves an $n$-dimensional plane $P'_1$ defined  over $\Q$, we apply an element of
$C_G(\Z)$ to find another $\rho(G)$-invariant $n$-plane $P'_2$ defined over $\Q$ 
which is in direct sum with $P'_1$. Thus, there is an element $B$ of $\SL_{2n}(\Q)$ 
which maps the standard decomposition $\R^{2n}=\R^n\oplus \R^n$ 
to $\R^{2n}= P'_1\oplus P'_2$, and conjugates $\rho(G)$  to the diagonal copy
of $\SL_n(\R)$ in $\SL_{2n}(\R)$. The group $\rho(\Gamma)$ is virtually conjugate, by the same matrix $B$, 
to the intersection of the diagonal copy of $\SL_n(\R)$ with $\SL_{2n}(\Z)$, so that $\Gamma$ 
is commensurable to $\SL_n(\Z)$ in $\SL_n(\R)$. This shows the following implications
\[
(1) \Rightarrow (2) \Rightarrow (3) \Rightarrow  (4) .
\]
Assume (4). 
Let $T$ be an element of $\SL_n(\R)$ such that $T\Gamma T^{-1}$ intersects
$\SL_n(\Z)$ on a finite index subgroup.
Let $N^+$ be the group of upper triangular matrices in $\SL_n(\R)$, let
$N^+(\Z)$ be its intersection with $T\Gamma T^{-1}$, and $N^+_\Gamma= T^{-1}N^+(\Z) D$. The action of $\rho(N^+_\Gamma)$ on $V$ fixes a $2$-plane $U$
pointwise, and $U$ is define over $\Z$ because $N^+_\Gamma$ preserves the lattice 
$\Lambda= \Z^{2n}$. Since $N^+_\Gamma$ fixes a unique direction $D$ in $P=\R^n$, 
the plane $U$ is equal to $D\times D$ in $V=P\times P$. Let $u=(a x_0, b x_0)$ be an element
of $U\cap \Lambda\setminus \{0\}$. Then the $n$-dimensional plane $P'$ defined by 
\[
P' = \left\{ (x,y)\in P\times P\, \vert \, ay=bx \right\}
\]
is $\rho(G)$-invariant and contains a lattice point $u\in \Lambda$. The orbit of $u$ under
$\rho(\Gamma)$ is a lattice in $P'$, and thus $P'$ is defined over $\Z$. This shows that 
\[
(4) \Rightarrow (1).
\]

Assume (2) and denote by $(x_1, \ldots, x_n,y_1, \ldots, y_n)$ the coordinates in 
$V=\R^n\oplus \R^n$ in which $\rho(G)$ is a diagonal copy of $\SL_n(\R)$. Then $C_G$
acts as $(x_i, y_i)\mapsto (ax_i+by_i, cx_i+by_i)$. Since this coordinates are
defined over $\Q$,  $C_G(\Z)$ is commensurable to $\SL_2(\Z)$, so that
(2) implies (5). Assume (5), 
and take a unipotent element $U$ in $C_G(\Z) \setminus \{{\text{Id}}\}$. The 
set of fixed points of $U$ in $V$ is an $n$-dimensional plane defined over $\Z$ which is
invariant by $\rho(G)$. Thus (5) implies (1), and all five properties are equivalent.
\end{proof}


\subsection{Stabilizers, cocompactness, and odd dimensions}

\subsubsection{} Let us now fix a non-zero element $(x_0,y_0)\neq 0$ in the intersection $\Lambda\cap (P \times P)$. 

\begin{rem}\label{rem:slnz}
If $x_0$ is proportional to $y_0$, with $bx_0=ay_0$, the $n$-plane $P'$ 
given by the equation $bx=ay$ is $G$-invariant and contains a
lattice point. As seen in the proof of Lemma \ref{lem:pepr}, this implies that $P'$ 
is defined over $\Q$ and that $\Gamma$ is commensurable to $\SL_n(\Z)$.
\end{rem}

We now assume that $x_0$ and $y_0$ are not collinear. Let $H$ be the stabilizer
of $(x_0,y_0)$ in $G$. Taking $x_0$ and $y_0$ as the first elements of a basis for $P$, 
the group $H$ can be identified with the semi-direct product $\SL_{n-2}(\R) \ltimes  \R^{2(n-2)}$
of matrices
\[
\left(
\begin{array}{ccc}
1 & 0 & u \\
0 & 1 & v \\
0 & 0 & A
\end{array}
\right)
\]
where $A$ is in $\SL_{n-2}(\R)$ and $u$ and $v$ are row vectors in $\R^{n-2}$. 
The $G$-orbit of $(x_0,y_0)$ in $V$ is homeomorphic to $G/H$ and its $\Gamma$-orbit 
is a discrete subset of $G/H$. From \cite{Margulis:1991}, Theorem 3.11 and Remark 3.12 (see also \cite{Shah:1991}, Lemma 2.8), we deduce that 
\[
\Gamma_H= \Gamma\cap H
\] is a lattice
in $H$. Since $H$ is the semi-direct product of its radical $N_H=\R^{2(n-2)}$ with the
semi-simple factor $S_H=\SL_{n-2}(\R)$, $\Gamma$ intersects $N_H$ onto a lattice $\Lambda_N$
(see \cite{Raghunathan:book}, Corollary 8.28 page 150). In particular, $\Lambda_N$ and $\Gamma$
contain unipotent elements; this implies that $\Gamma$ is not cocompact,
and proves the following proposition.

\begin{pro}
Let $\Gamma$ be a lattice in $G=\SL_n(\R)$, with $n\geq 3$. If $\Gamma$ preserves a lattice $\Lambda$
in the diagonal representation of $G$ in $\R^n\times \R^n$ then $\Gamma$ is not cocompact.
\end{pro}

\subsubsection{} Let us now denote by $E_0$ the $2$-plane contained in $P$ which is generated by $x_0$ and $y_0$. 
By construction, the $4$ dimensional space
\[
B_0=E_0\times E_0 \subset P\times P
\] 
coincides with the set of fixed points of $H$. Since
$\Gamma$ intersects $H$ on a lattice, $E_0\times E_0$ is defined over $\Q$.
Note that $B_0$  is $J$-invariant, i.e. is a complex subspace of complex dimension $2$ in $V$.

If we reproduce the same construction for another point $(x_1,y_1)$ in $P\times P \cap \Lambda$ 
with $x_1$ and $y_1$ not collinear, we get another $2$-plane $E_1$ in $P$ and another $4$-plane $B_1$
in $P\times P$. Let $k$ be the maximum number of such $4$-planes
$B_i= E_i\times E_i$, $0\leq i \leq k-1$, such that the sum of the $B_i$ has dimension $4k$. Let ${\mathbf{B}}_c$ 
be the direct sum of the $B_i$; this plane of dimension $4k$ is defined over $\Q$. Similarly, the direct
sum ${\mathbf{E}}_c $ of the $E_i$ has dimension $2k$.

Let $(x_{k}, y_k)\neq (0,0)$ be an element of $\Lambda$ such that $x_k$ and $y_k$
are not proportional. 
The intersection of the corresponding $4$-plane $B_k$ with ${\mathbf{B}}_c $ has positive dimension, and $E_k$
intersects ${\mathbf{E}}_c $. If this intersection is a line $D$, we see that $B_k\cap {\mathbf{B}}_c $ is
equal to $D\times D$; since  $B_k\cap {\mathbf{B}}_c $ is rational, there exists a point $(u,v)\neq 0$  in $\Lambda\cap(D\times D)$. 
The vectors $u$ and $v$ are proportional, and Remark \ref{rem:slnz} implies that $\Gamma$
is commensurable to $\SL_n(\Z)$.
Thus, we can assume that all planes $B_k$, for all starting points $(x_k,y_k)$ in $\Lambda$ are 
indeed contained in ${\mathbf{B}}_c $. This shows the following. 

\begin{lem}
Let $\Gamma$ be a lattice in $\SL_n(\R)$ if $\Gamma$ preserves a lattice in 
the diagonal representation of $\SL_n(\R)$ in $\R^{2n}$, then 
\begin{itemize}
\item either $\Gamma$ is commensurable to $\SL_n(\Z)$
\item or $n=2k$ and there exists $k$ distinct points $(x_j,y_j)$ in $\Lambda$ 
such that the $4$-planes $B_j$ constructed above are in direct sum.
\end{itemize}
\end{lem}

\begin{rem}
In particular, if the dimension $n$ is odd, the lattice $\Gamma$ is commensurable to 
$\SL_n(\Z)$ and the group $G$ is conjugate by an element of $\SL_{2n}(\Q)$ to the diagonal
copy of $\SL_n(\R)$ in $\SL_{2n}(\R)$.
\end{rem}

\subsection{Quaternion algebras and even dimensions}\label{par:quaternions}

We now explain, conserving the same notation, how all examples can be constructed
in even dimension. 

\subsubsection{Quaternion algebras and lattices in $\SL_2(\R)$ (see \cite{Katok:book}, \cite{Witte:prebook})} \label{par:quater}

Let $a$ and $b$ be two integers. Let $\H_{a,b}$ (or $\H_{a,b}(\Q)$) be the quaternion algebra over the rational 
numbers $\Q$ defined by its basis $(1,\ii,\jj,\kk)$, with 
\[
\ii^2=a, \, \jj^2=b, \, \ii\jj=\kk=-\jj\ii.
\]
This algebra embeds into the space of $2\times 2$ matrices over $\Q(\sqrt{a})$ by mapping
$\ii$ and $\jj$ to the matrices
\[
\left(
\begin{array}{cc}
\sqrt{a} & 0 \\ 
0 & -\sqrt{a}
\end{array}
\right), \quad
\left(
\begin{array}{cc}
0 & 1 \\ 
b & 0
\end{array}
\right).
\]
In what follows, we denote by $\H_{a,b}(\Z)$ the set of quaternions with coefficients in $\Z$, 
and by $\H_{a,b}(\R)$ the tensor product $\H_{a,b}\otimes_\Q \R$.
The determinant of the matrix which is associated to a quaternion $q=x + y \ii + z \jj + t\kk$ 
is   equal to its reduced norm 
\[
\Nrd(q) = x^2 - a y^2 - b z^2 + ab t^2.
\]
Assume that $\H_{a,b}$ is a division algebra, i.e. that $\Nrd(q)\neq 0$ if $q\neq 0$ is an element
of $\H_{a,b}(\Q)$. Then the image of $\H_{a,b}(\Q)^*$ is contained in $\GL_2(\Q(\sqrt{a}))$; moreover
\begin{enumerate}
\item[(1)] The group of quaternions $q$ with reduced norm $1$ and integer coefficients determines
a cocompact lattice $C_{a,b}$ in $\SL_2(\R)$;
\item[(2)] This lattice acts by left multiplication on $\Mat_2(\R)\simeq \R^4$, preserving the (image of the) lattice $\H_{a,b}(\Z)$.
\end{enumerate}
Quaternions also act by right multiplication. The group of invertible linear transformations of the vector
space $\H_{a,b}(\R) =  \Mat_2(\R)$
that commute with the left action of $C_{a,b}$ coincides with the group $\H_{a,b}(\R)^*$ of quaternions with 
real coefficients and non-zero reduced norms, that is with the group $\GL_2(\R)$, acting on $\Mat_2(\R)$ 
by right multiplications. The third property we need is the following.
\begin{enumerate}
\item[(3)] If $L\subset\H_{a,b}(\R)$ is a lattice which is invariant by a finite index subgroup of $C_{a,b}$,
then $L$ is commensurable to a right translate   $\H_{a,b}(\Z)B$ by an element $B$ of the centralizer $\GL_2(\R)$ of
$C_{a,b}$.
\end{enumerate}
The last important fact  characterizes the lattices $C_{a,b}$.
\begin{enumerate}
\item[(4)] If $C$ is an arithmetic lattice in $\SL_2(\R)$ with rational  traces ${\text{tr}}(c)\in \Q$ for all $c\in C$, 
then $C$ is commensurable to $C_{a,b}$ for some division algebra $\H_{a,b}$ with $a$ and $b$ in $\Z$.
\end{enumerate}

\subsubsection{} Let us come back to the study of $C_G(\Z)$ and its action on $V$ and the $B_j$.
Since $C_G(\Z)$ preserves $\Lambda$ and the representation $C_G(\R)\to \GL(V)$ is equivalent
to $n$ diagonal copies of the standard   action of $\SL_2(\R)$ on $\R^2$, one concludes that $C_G(\Z)$ 
is an arithmetic lattice with rational traces. Thus, property (4) in Section \ref{par:quater} implies that $C_G(\Z)$
is commensurable to $C_{a,b}$ for some division algebra $\H_{a,b}$. Moreover, property (3) shows that 
$\Lambda_j$ is commensurable to a left translate of $\H_{a,b}(\Z)$ by an element of the centralizer of $C_G(\Z)$. 
This shows that the lattice $\Lambda$ itself is commensurable to $\H_{a,b}(\Z)^{k}$ up to the action of an element
in the group of complex linear transformations of $V$. 

The group $\rho(\Gamma)$ acts on $V$, preserves $\Lambda$, and commute to $C_G$. Thus, up to a linear isomorphism, 
$\rho(\Gamma)$ is commensurable to the group of linear transformation $L_k({a,b};\Z)$ of $\H_{a,b}(\R)^{k}$ which
preserves the lattice $\H_{a,b}(\Z)^{k}$ and commutes to the diagonal action of $\H_{a,b}(\R)^*$ by right multiplications on $\H_{a,b}(\R)^{k}$.

\begin{thm}
Let $M$ be a complex torus of complex dimension $n$, together with a faithful holomorphic action of a lattice $\Gamma$ of $\SL_n(\R)$.
Then
\begin{itemize}
\item either $M$ is isogeneous to the product of $n$ copies of an elliptic curve $\C/\Lambda_0$ and $\Gamma$ is commensurable
to $\SL_n(\Z)$;
\item or $n=2k$ is even and there exists a division algebra  $\H_{a,b}$ over $\Q$ such that  $M$ is isogeneous to 
the product of $k$ copies of the abelian surface $\C^2/\H_{a,b}(\Z)$ and $\Gamma$ is commensurable to the group of 
automorphisms of the abelian group $\H_{a,b}(\Z)^k$ that commute to the diagonal action of $\H_{a,b}(\Z)$ by left multiplications.\end{itemize}
In particular, $\Gamma$ is not cocompact and $M$ is an abelian variety.
\end{thm}

\subsection{Kummer examples and singularities}\label{part6:Kummer}
 Once we have the list of possible tori and lattices, Kummer examples
are obtained by a quotient with respect to a finite group of automorphisms of the torus. 

Let $A=\C^n/\Lambda$ be a torus and $\Gamma$ be a lattice in $\SL_n(\R)$ or $\SL_n(\C)$ acting faithfully on $M$.
Let $F$ be a finite group of automorphisms of $A$ which is normalized by the action of $\Gamma$. From 
Lemma \ref{lem:Fcyclic}, we can assume that $F$ is a finite cyclic group of homotheties.

If $A$ is isogeneous to $(\C/\Lambda_0)^n$, with $\Lambda_0$ a lattice in $\C$ and $\Lambda=\Lambda_0^n$, the
order of $F$ is $1$, $2$, $3$, $4$ or $6$ (see \cite{Cantat-Zeghib}). If $n=2k$ is even and $M$ is isogeneous to $(\C^2/\H_{a,b}(\Z))^k$,
the same conclusion holds: The finite group $F$ is contained in the centralizer of $\Gamma$, that is in the 
group $C_G$, preserves $\Lambda$, and is finite cyclic. Thus, $F$ can be identified to a finite cyclic subgroup
of $\C_{a,b}$. Viewed as a subgroup of $\SL_2(\R)$, the traces are even integers, and thus finite order elements
have trace in $\{-2, 0,  2\}$. Thus  the order of the cyclic group $F$ is bounded by $2$ in this case. 

This proves the following fact.

\begin{pro}
Let $M_0$ be Kummer orbifold $A/F$ where $A=\C^n/\Lambda$ is a torus of dimension $n$ and
$F$ is a finite group of automorphisms of $A$. Assume that there is a faithful action of a lattice
in an almost  simple Lie group $G$ of a rank $n-1$ on $M_0$. 
Then $M_0$ is the quotient $A'/F'$ of a torus $A'$ isogenous to $A$ by a finite cyclic group 
$F'$ which is generated by a 
scalar multiplication 
\[
(x_1, \ldots , x_n)\mapsto (\eta x_1, \ldots, \eta x_n)
\]
where $\eta$ is  a root of unity of order $2$, $3$, $4$ or $6$.  
\end{pro}

%
%

 

\end{document}